\documentclass[margin=1cm]{article}

\usepackage{amsmath,amssymb}
\usepackage{graphics} 
\usepackage{epsfig}
\usepackage{graphicx}  
\usepackage{subfigure}
\usepackage{epstopdf} 
\usepackage{placeins}
\usepackage[colorlinks=true]{hyperref}
\usepackage{fancyhdr}

\hypersetup{urlcolor=blue, citecolor=red}

\usepackage{theorem}
\newtheorem{proposition}{Proposition}

\usepackage{changes}
\usepackage{todonotes}

\renewcommand{\emph}[1]{{\it #1}}

\newcommand{\Pe}{\text{Pe}}
\def\VhE{V_h(E)}
\def\VhO{V_h(\Omega_h)}

\newcommand{\numberset}{\mathbb}
\newcommand{\N}{\numberset{N}}
\newcommand{\R}{\numberset{R}}
\newcommand{\Pk}{\numberset{P}}

\newcommand{\nn}{\boldsymbol{n}}

\newcommand{\xx}{\boldsymbol{x}}
\newcommand{\uu}{\boldsymbol{u}}

\usepackage{authblk}
\usepackage{setspace}

\begin{document}

\title{$C^1$-VEM for some variants of the Cahn-Hilliard equation: a numerical exploration} 
\author[1]{Paola F. Antonietti\thanks{paola.antonietti@polimi.it}}
\author[2]{Simone Scacchi\thanks{simone.scacchi@unimi.it}}
\author[3]{Giuseppe Vacca\thanks{giuseppe.vacca@uniba.it}}
\author[1]{Marco Verani\thanks{marco.verani@polimi.it}}
\affil[1]{MOX-Dipartimento di Matematica, Politecnico di Milano,\newline Piazza Leonardo da Vinci 32, 20133 Milano, Italy}
\affil[2]{Dipartimento di Matematica, Universit\`a degli Studi di Milano, \newline Via Saldini 50, 20133 Milano, Italy}
\affil[3]{Dipartimento di Matematica, Universit\`a degli Studi di Bari, \newline Via Edoardo Orabona 4, 70125, Italy}

\maketitle

\begin{abstract}
We consider the $C^1$-Virtual Element Method (VEM) for the conforming numerical approximation of some variants of the Cahn-Hilliard equation on polygonal meshes. In particular, we focus on the discretization of the advective Cahn-Hilliard problem and the Cahn-Hilliard inpainting problem. We present the numerical approximation and several numerical results to assess the efficacy of the proposed methodology.
\end{abstract}

\vspace{0,5cm}


\section{Introduction}\label{sec:intro}
The Cahn-Hilliard equation, which is a fourth-order nonlinear parabolic problem, was  initially  introduced  as  a  diffusive  interface  model  to  characterize  the  phase  segregation of binary alloys at constant temperature \cite{Cahn-Hilliard:1958}. Compared to sharp-interface models, where the individual interfaces need to be explicitly tracked,  the advantage of a diffuse-interface approach is that topological changes are automatically handled, since interfaces are treated in a diffuse manner thanks to the introduction of a parameter which, varying continuously, accounts for the different material phases and/or the concentration of the different components.  
Since the seminal paper by Cahn and Hilliard, several different variants have been studied (see, e.g., the book \cite{Miranville:book} and the references therein) covering  a wide spectrum of applications. Here we mention, for example, the modeling and simulation of solidification processes, spinodal decomposition, coarsening of precipitate phases, shape memory effects, re-crystallisation, dislocation dynamics \cite{Chen2002113,Emmerich20081,Moelans2008268,Steinbach2009}, 
wettability~\cite{Frank2018282}, diblock copolymer~\cite{Thomas1988598,Regazzoni2018562}, 
tumor growth~\cite{Agosti20177598,Frigeri2015215,Hawkins-Daarud20123,Wise2008524,Chatelain2011}, 
image inpainting~\cite{Bertozzi2007285,Bertozzi2007913}, 
crystal growth~\cite{Elder20022457011,Emmerich20111,Wang2011945,
Grasselli20161523} and crack propagation~\cite{Miehe20101273,Kuhn20103625,
Borden201277}.\\
    
In the last decades, different numerical techniques have been utilized to solve the Cahn-Hilliard equation and its variants, including finite difference, finite element, and spectral methods. A crucial difficulty in designing numerical schemes is that these equations typically involve spatial differential operators that are higher than second-order. Therefore, standard conforming $C^0$  Finite Element Methods (FEMs) are ruled out and approximation spaces with higher regularity are required.  However, the construction of such approximation spaces with higher regularity is deemed a difficult task because it requires a set of basis functions with such a global regularity.
Examples in this direction can be found all along the history of
finite elements: from the oldest works in the sixties of the last
century,
e.g.,~\cite{Argyris-Fried-Scharpf:1968,Bell:1969,Clough-Tocher:1965}
to  the most recent attempts
in~\cite{Zhang:2009,Zhang:2016,Zhang-Hu:2015,Wu-Lin-Hu:2021}.
Despite its intrinsic difficulty, designing approximations with global
$C^{1}$ or higher regularity is still a major research topic. In the literature there is a 
limited number of works addressing  the solution of the Cahn-Hilliard equation by the $C^1$-FEM, see~\cite{Elliott1986339,Elliott198797}.
To circumvent the  well known difficulty met in the implementation of $C^1$-FEMs, another possibility
is the use of non-conforming (see, e.g., \cite{Elliott-French:1989}) or discontinuous (see, e.g., \cite{Wells2006860,Engel20023669})
methods; obviously in such cases the discrete solution will not satisfy $C^1$ regularity. \\

Alternatively, the Cahh-Hilliard problem can be split into a coupled of lower-order differential equations and mixed formulation can be employed for discretization at the expense of introducing additional unknowns, see, e.g., \cite{Feng200447,Kay-Styles-Suli:2009,Aristotelous20132211,Liu20191509,Chave20161873, Liu2020}. 
Recently, in \cite{Gomez20084333} isogeometric analysis has been employed to discretize the Cahn-Hilliard problem, whereas  in \cite{MR3062036} the same approach has been used to discretize the advective Cahn-Hilliard equation. A remarkable feature of this methodology is that the approximation spaces exhibit higher-order continuity properties, thus avoiding the use of mixed formulations.  More recently, in \cite{Antonietti-BeiraodaVeiga-Scacchi-Verani:2016} the $C^1$-Virtual Element Method (VEM) has been employed to discretize the Cahn-Hilliard equation on polygonal meshes, employing highly regular conforming approximation spaces, thus circumventing the introduction of additional variables typical of mixed formulations.\\

Roughly speaking, the VEM is a Galerkin-type projection method that
generalize the finite element method, which was originally designed
for simplicial and quadrilateral/hexahedral meshes, to polygonal/polyhedral (polytopal, for short)
meshes. The VEM has been originally introduced in~\cite{BeiraodaVeiga-Brezzi-Cangiani-Manzini-Marini-Russo:2013} and 
does not require the explicit knowledge of the basis functions
spanning the approximation space. The functions that belong to such approximation spaces are dubbed as
``virtual'' as they are never really computed, with the noteworthy
exception of a subspace of polynomials that are indeed used in the
formulation and implementation of the method.
The virtual element functions are uniquely characterized by a set of
values, the so called \emph{degrees of freedom}. The VEM can then be implemented 
using only the degrees of freedom and the polynomial part of the approximation space.
The crucial idea behind the VEM is that the elemental approximation
space is defined elementwise as the solution of a partial
differential equation. Then the global approximation space is obtained by globally ``\emph{gluing}'' the local spaces in an arbitrary highly regular conforming way. 
Thus, the virtual element ``\emph{paradigm}'' provides a major breakthrough as it allows to obtain highly-regular Galerkin methods, and  the construction of numerical approximation sof any order of accuracy  on unstructured two-dimensional and three-dimensional meshes  made by general polytopal elements.

The first works proposing a $C^{1}$-regular conforming VEM addressed the
classical plate bending
problems~\cite{Brezzi-Marini:2013,Chinosi-Marini:2016}, second-order
elliptic
problems~\cite{BeiraodaVeiga-Manzini:2014,BeiraodaVeiga-Manzini:2015},
and the nonlinear Cahn-Hilliard
equation~\cite{Antonietti-BeiraodaVeiga-Scacchi-Verani:2016}.
More recently, highly regular virtual element spaces were considered
for the von {K}\'arm\'an equation modeling the deformation of very
thin plates~\cite{Lovadina-Mora-Velasquez:2019}, geostrophic
equations~\cite{Mora-Silgado:2021} and fourth-order subdiffusion
equations~\cite{subdiffusion:2021}, two-dimensional plate vibration
problems of Kirchhoff plates~\cite{Mora-Rivera-Velasquez:2018}, 
transmission eigenvalue problems~\cite{Mora-Velasquez:2018}, and
fourth-order plate buckling eigenvalue
problems~\cite{Mora-Velasquez:2020}.
In~\cite{Antonietti-Manzini-Verani:2019}  the
highly-regular conforming VEM for the two-dimensional polyharmonic
problem $(-\Delta)^{p_1} u= f $, $p_1\geq 1$ has been proposed.
The VEM is based on an approximation space that locally contains
polynomials of degree $r\geq 2p_1-1$ and has a global $H^{p_1}$
regularity.
In~\cite{Antonietti-Manzini-Scacchi-Verani:2021},  this
formulation has been extended to a virtual element space that can have arbitrary
regularity $p_2\geq p_1\geq1$ and contains polynomials of degree
$r\geq p_2$.
VEMs for three-dimensional problems are also available for the
fourth-order linear elliptic equations~\cite{BeiraodaVeiga-Dassi-Russo:2020} (see
also~\cite{Brenner-Sung:2019}), and highly-regular conforming VEM in any dimension 
has been proposed in~\cite{Huang:2021}.\\

In this paper, hinging upon the use of $C^1$-VEM,  we study the conforming virtual element approximation on polygonal meshes of two variants of the Cahn-Hilliard equation, namely the Advective Cahn-Hilliard (ACH) problem and the Cahn-Hilliard Inpainting problem (CHI). Those variants have been selected both for their relevance in applications and for the presence, with respect to the classical Cahn-Hilliard equation,  of the additional convective term in the ACH problem and the reaction term in the CHI problem.  The numerical treatment of those terms is new in the context of the conforming virtual element discretization of Cahn-Hilliard equations.
It is also worth mentioning that the numerical treatment of the advective Cahn-Hilliard represents an important  preliminary step to tackle in future works the virtual element approximation of  more complicated  problems, as the convective nonlocal Cahn-Hilliard (see, e.g.,~\cite{DellaPorta20151529}) or the Navier-Stokes-Cahn-Hilliard problem (see, e.g., \cite{Gal2010401,Gal2016} and \cite{Bao20128083,Feng20061049,Kay200815}).

\medskip

The paper is organized as follows.  In Section~\ref{sec:2} we introduce the continuous problems, whereas in Section~\ref{sec:3} we present their conforming virtual element approximation. In Section~\ref{sec:num} we collect and discuss several numerical results to show the efficacy of our discretization methodology. Finally, in Section~\ref{S:conclusions} we draw some conclusions.\\

\noindent {\bf Notation}. Throughout the paper, we will follow the usual notation for Sobolev spaces
and norms \cite{Adams:1975}.
Hence, for an open bounded domain $\omega$,
the norms in the spaces $W^s_p(\omega)$ and $L^p(\omega)$ are denoted by
$\|{\cdot}\|_{W^s_p(\omega)}$ and $\|{\cdot}\|_{L^p(\omega)}$, respectively.
The norm and seminorm in $H^{s}(\omega)$, $s\geq 1$, are denoted  by
$\|{\cdot}\|_{s,\omega}$ and $|{\cdot}|_{s,\omega}$, respectively.
The $L^2$-inner product and the $L^2$-norm are denoted by $(\cdot,\cdot)_{\omega}$ and $\|\cdot\|_{\omega}$, respectively. The subscript $\omega$ may be omitted when $\omega$ is the whole computational
domain $\Omega$.
We denote with $\xx = (x_1, \, x_2)$ the independent variable. 
With the usual notation the symbols $\nabla$, $\Delta$, $\Delta^2$, $D^2$ 
denote the gradient, the laplacian, the bilaplacian and the Hessian for (regular enough) scalar functions, whereas $\partial_t$ denotes the derivative with respect to the time variable.

%

\section{Continuous problems}
\label{sec:2}
In this Section~we introduce the two variants of the classical Cahn-Hilliard problem, whose numerical discretization  will be addressed in the sequel of this paper. More specifically, we consider the Advective Cahn-Hilliard  problem and the Cahn-Hilliard Impainting problem. For each variant, we provide the weak formulation that will the basis for the construction of the virtual element discretization. 
 
Let $\Omega\subset \mathbb{R}^2$ be an open bounded domain. Let  
$\psi:\mathbb{R}\to\mathbb{R}$ with $\psi(x)= (1-x^2)^2/4$ and let $\phi(x) = \psi^\prime (x)$,  we consider the following two variants of the Cahn-Hilliard problem, where $\gamma \in {\mathbb R}^+$, $0 < \gamma \ll 1$, represents the interface parameter.
\\

\noindent {\bf Advective Cahn-Hilliard problem}.
For a given final time $T>0$, find $c(x,t):\Omega \times [0,T] \rightarrow {\mathbb R}$ such that:
\begin{equation}\label{eq:CHforte:1}
\left\{
\begin{aligned}
& \partial_t c - \frac{1}{\Pe}\Delta \big( \phi(c) - \gamma^2 \Delta c \big) + {\rm div} (\uu c) = 0 
&& \textrm{ in } \Omega \times (0,T], \\
& c(\cdot, 0) = c_0(\cdot) 
&& \textrm{ in } \Omega, \\
& \partial_{\bf n} c = {  \partial_{\bf n}} \big( \phi(c) - \gamma^2 \Delta c \big) = 0
&& \textrm{ on } \partial\Omega \times (0, T],
\end{aligned}
\right.
\end{equation}
where $\partial_{\nn}$ denotes the (outward) normal derivative and $\Pe$ is a positive constant. {We note that on the boundary of the domain we impose no-flux type condition both on $c$ and on the so-called chemical potential $\phi(c) - \gamma^2 \Delta c $}. 
Finally, $\uu \in H({\rm div}, \Omega) \cap [C^0(\Omega)]^2$ is a given function such that ${\rm div} \uu=0$ in $\Omega$ and $\uu \cdot \nn=0$ on $\partial \Omega$. Here 
$$H({\rm div}, \Omega)=\{\boldsymbol{v}\in [L^2(\Omega)]^2 : {\rm div} \boldsymbol{v} \in   L^2(\Omega)\}\,.$$

\noindent {\bf Cahn-Hilliard inpainting problem}. Let $f$ be a given binary image and $D\subset\Omega$ be the inpainting domain. For a given final time $T>0$, find $c(x,t):\Omega \times [0,T] \rightarrow {\mathbb R}$ such that:
\begin{equation}\label{eq:CHforte:2}
\left\{
\begin{aligned}
& \partial_t c - \Delta \big( \frac{1}{\gamma}\phi(c) - \gamma \Delta c \big) + \lambda(x)(f-c) = 0 
&& \textrm{ in } \Omega \times (0,T], \\
& c(\cdot, 0) = c_0(\cdot) 
&& \textrm{ in } \Omega, \\
&  \partial_{\bf n} c = {  \partial_{\bf n}} \big( \frac 1 \gamma\phi(c) - \gamma \Delta c \big) = 0
&& \textrm{ on } \partial\Omega \times (0, T],
\end{aligned}
\right.
\end{equation}
where   
 $$
\lambda(x)=
\left\{
\begin{aligned}
&\lambda_0 ,\quad &\xx\in \Omega\setminus D,\\
&0, \quad & \xx \in D,
\end{aligned}
\right.
$$
$\lambda_0$ being a positive parameter. See, e.g., \cite{Bertozzi2007285,Bertozzi2007913} for more details on the model. 

We now briefly introduce the variational formulations  of \eqref{eq:CHforte:1} and \eqref{eq:CHforte:2} that will be used to derive the virtual element discretizations. To this aim, we preliminary define the following bilinear forms 
\begin{equation}
\label{eq:forms}
\begin{aligned}
& a^{D^2}(v,w) = \int_\Omega (D^2 v) : (D^2 w) \, {\rm d}\Omega
&& \forall v,w \in H^2(\Omega), \\
& a^0(v,w) = \int_\Omega v \, w \, {\rm d}\Omega 
&& \forall v,w \in L^2(\Omega),\\
& b(v,w) = \int_\Omega {\bf u}  \cdot \nabla v \ w \, {\rm d}\Omega 
&&\forall v,w \in H^1(\Omega),\\
\end{aligned}
\end{equation}
and the semi-linear forms 
\begin{equation}
\label{eq:formar}
\begin{aligned}
&l(f;v,w) = \int_\Omega \lambda (f - v)w \, {\rm d}\Omega
&&\forall v,w \in L^2(\Omega) \,,
\\
&r(z;v,w) = \int_\Omega \phi'(z) \nabla v \cdot \nabla w \, {\rm d}\Omega
&&\forall z,v,w \in H^2(\Omega) \,.
\end{aligned}
\end{equation}
Finally, we introduce  the space
\begin{equation}\label{V}
V = \big\{ v \in H^2(\Omega) \: : \: \partial_{\bf n} { v} = 0 \textrm{ on } \partial \Omega \big\}.
\end{equation}
The weak formulation of problem \eqref{eq:CHforte:1} reads as follows: find $c(\cdot, t) \in V$ s.t.
\begin{equation}\label{contpbl:1}
\left\{
\begin{aligned}
& a^0(\partial_t c,v) + \frac{ \gamma^2}{\Pe} a^{D^2}(c,v) + \frac{1}{\Pe}r(c;c,v) + b(c,v)= 0 \quad \forall v \in V ,   \\
& c(\cdot,0)=c_0  \,.
\end{aligned}
\right.
\end{equation}
Similarly, the weak formulation of problem \eqref{eq:CHforte:2} reads as follows: find $c(\cdot, t) \in V$ s.t.
\begin{equation}\label{contpbl:2}
\left\{
\begin{aligned}
& a^0(\partial_t c,v) + { \gamma} a^{D^2}(c,v) + \frac{1}{\gamma} r(c;c,v) + l(f; c, v)= 0 \quad \forall v \in V , 
\\
& c(\cdot,0)=c_0 \,.
\end{aligned}
\right.
\end{equation}

\section{Virtual element discretization}
\label{sec:3}
In this Section~we describe the virtual element discretization of problems \eqref{contpbl:1}-\eqref{contpbl:2} on computational meshes made of  general polygons. In particular, in Section~\ref{sub:mesh} we introduce the assumptions on the regularity of the polygonal mesh together with the definition of crucial projector operators that will be fundamental in the construction of the virtual element discretization. In Section~\ref{sub:spaces} we describe the $C^1$-Virtual Element spaces that will of paramount importance to guarantee a conforming approximation of the Cahn-Hilliard problems. Finally, in Section~\ref{sub:vem problem} we introduce the semi-discrete in space virtual element discretization of  \eqref{contpbl:1}-\eqref{contpbl:2} together with a fully discrete scheme based on the use of the backward Euler method for time discretization.
\subsection{Mesh assumptions and polynomial projections}
\label{sub:mesh}

From now on, we will denote with $E$ a general polygon,  having $n_e$ edges $e$, moreover $|E|$ and $h_E$ will denote the area  and the diameter of $E$, respectively.
Let $\{\Omega_h\}_h$ be a sequence of decompositions of $\Omega$ into general polygons $E$, 
where the granularity $h$ is defined as $h = \sup_{E \in \Omega_h} h_E$.
We suppose that $\{\Omega_h\}_h$ fulfills the following assumption:\\
\textbf{(A1) Mesh assumption.}
There exists a positive constant $\rho$ such that for any $E \in \{\Omega_h\}_h$   
\begin{itemize}
\item Any $E \in \{\Omega_h\}_h$ is star-shaped with respect to a ball $B_E$ of radius $ \geq\, \rho \, h_E$;
\item Any edge $e$ of any $E \in \{\Omega_h\}_h$ has length  $ \geq\, \rho \, h_E$.
\end{itemize}
We remark that the hypotheses above, though not too restrictive in many practical cases, could possibly be further relaxed, combining the present analysis with the studies in~\cite{BLR:2017,brenner-guan-sung:2017,brenner-sung:2018}.\\

Referring to Problem \eqref{contpbl:2}, we assume that for any $h$ there exists $D_h \subseteq \Omega_h$ such that $D_h$ is a decomposition of $D$, i.e. $\Omega_h$ matches with the subdivision of $\Omega$ into $D$ and $\Omega \setminus D$. \\

We denote with $\Sigma_h$ the set of all the mesh edges and for any $E \in \Omega_h$ we denote with $\Sigma_h^E$ the set of the edges of $E$.
Furthermore for any mesh vertex $\boldsymbol{\xi}$ we denote with $h_{\boldsymbol{\xi}}$ the average of the diameters of the elements having $\boldsymbol{\xi}$ as a vertex.
The total number of vertexes, edges and elements in the decomposition $\Omega_h$ are denoted  by $N_V$, $N_e$ and $N_P$, respectively.\\

Using standard VEM notations, for any mesh object $\omega \in \Omega_h \cup \Sigma_h$ and for any
$n \in \N$  let us introduce the space $\Pk_n(\omega)$ to be the space of polynomials defined on $\omega$ of degree $\leq n$  (with the extended notation $\Pk_{m}(\omega)=\{ 0 \}$ for any negative integer $m$). Moreover, 
$\widehat{\Pk}_{n\setminus m}(\omega)= \Pk_n(\omega) \setminus \Pk_m(\omega)$, for $m <n$, denotes the polynomials in $\Pk_n(\omega)$ with monomials of degree strictly greater than $m$. Finally, we introduce 
 the broken polynomial space
\[
\Pk_n(\Omega_h) = \{q \in L^2(\Omega) \quad \text{s.t.} \quad q|_E \in  \Pk_n(E) \quad \text{for all $E \in \Omega_h$}\} \,.
\]
For any non-negative $s\in \R$ let us introduce the broken space:
\[
H^s(\Omega_h) = \{v \in L^2(\Omega) \quad \text{s.t.} \quad v|_E \in  H^s(E) \quad \text{for all $E \in \Omega_h$}\} \,.
\]
Furthermore, we introduce the following notation: let $\{\mathcal{X}^E\}_{E \in \Omega_h}$ be a family 
of forms 
$\mathcal{X}^E \colon \prod_{j=1}^\ell H^{s_j}(E) \to \R$,
then we define 
\begin{equation}
\label{eq:XG}
\mathcal{X}\colon \prod_{j=1}^\ell H^{s_j}(\Omega_h) \to \R \,,
\quad 
\mathcal{X}(u_1, \dots, u_\ell) = 
\sum_{E \in \Omega_h}\mathcal{X}^E(u_1, \dots, u_\ell) \,,
\end{equation}
for any $u_j \in H^{s_j}(\Omega_h)$, and $j=1, \dots, \ell$. 

\noindent
For any $E \in \Omega_h$,  let us introduce the following polynomial projections:
\begin{itemize}
\item the $\boldsymbol{L^2}$\textbf{-projection} $\Pi_n^{0, E} \colon L^2(E) \to \Pk_n(E)$, given by
\begin{equation}
\label{eq:P0_k^E}
\int_{E} q_n (v - \, {\Pi}_{n}^{0, E}  v) \, {\rm d} E = 0 \qquad  \text{for all $v \in L^2(E)$  and $q_n \in \Pk_n(E)$,} 
\end{equation} 
with obvious extension for vector functions $\Pi^{0, E}_{n} \colon [L^2(E)]^2 \to [\Pk_n(E)]^2$ and 
tensor functions $\Pi^{0, E}_{n} \colon [L^2(E)]^{2 \times 2} \to [\Pk_n(E)]^{2 \times 2}$;

\item the $\boldsymbol{H^2}$\textbf{-seminorm projection} ${\Pi}_{n}^{D^2,E} \colon H^2(E) \to \Pk_n(E)$, defined by 
\begin{equation}
\label{eq:PD2_k^E}
\left\{
\begin{aligned}
& \int_{E} D^2 q_n : D^2 ( v - \, {\Pi}_{n}^{D^2,E}   v)\, {\rm d} E = 0 \quad  \text{for all $v \in H^2(E)$ and  $q_n \in \Pk_n(E)$,} \\
& \int_{\partial E}(v - \,  {\Pi}_{n}^{D^2, E}  v) \, {\rm d}s= 0 \, ,
\\
& \int_{\partial E}\partial_{\nn}(v - \,  {\Pi}_{n}^{D^2, E}  v) \, {\rm d}s= 0 \, .
\end{aligned}
\right.
\end{equation}
\end{itemize}
The global counterparts of the previous projections
\[
\Pi_n^{0} \colon L^2(\Omega_h) \to \Pk_n(\Omega_h)\,, \quad 
{\Pi}_{n}^{D^2} \colon H^2(\Omega_h) \to \Pk_n(\Omega_h)
\]
are defined for all $E \in \Omega_h$ by
\begin{equation}
\label{eq:proj-global}
(\Pi_n^{0} v)|_E = \Pi_n^{0,E} v \,,
\qquad
(\Pi_n^{D^2} v)|_E = \Pi_n^{D^2,E} v \,.
\end{equation}

In the following the symbol $\lesssim$ will denote a bound up to a generic positive constant,
independent of the mesh size $h$, but which may depend on 
$\Omega$, on the ``polynomial'' order of the
method $k$ and on the regularity constant appearing in the mesh assumption \textbf{(A1)}.
\subsection{Virtual Element space}
\label{sub:spaces}

In the present Section~we outline an overview of the $H^2$-conforming Virtual Element space \cite{BM:2012,Antonietti-BeiraodaVeiga-Scacchi-Verani:2016,Antonietti-Manzini-Scacchi-Verani:2021} combined with the construction proposed in \cite{projectors} in order to define the 
``enhanced'' version of such space such that the ``full'' $L^2$-projection $\Pi_{k}^{0,E}$ is computable by the degrees of freedom (DoFs).

Let $k \geq 2$ be the  ``polynomial'' order of the method.
We thus consider on each polyhedral element $E \in \Omega_h$ the ``enhanced'' virtual space
\begin{equation}
\label{eq:v-loc}
\begin{aligned}
\VhE = \biggl\{  
v \in C^1(\overline{E}) \,\,\, \text{s.t.} \,\,\,
(i) & \, \,  \Delta^2    v  \in \Pk_{k}(E) \,, 
\\
(ii)&  
\,  \, v_{|e} \in \Pk_{\widetilde{k}}(e) \,\,\, \forall e \in \Sigma^E_h \,, 
\\
(iii) &
\, \partial_{\nn} v_{|e} \in \Pk_{k-1}(e) \,\,\, \forall e \in \Sigma^E_h, 
\\
(iv) &
\,  \int_E \bigl( v - \Pi^{D^2,E}_k v) \widehat{p}_{k} \, {\rm d}E = 0
\,\,\, \text{$\forall \widehat{p}_{k} \in \widehat{\Pk}_{k \setminus k-4}(E)$}
\,\,\biggr\} \,,
\end{aligned}
\end{equation}
where $\widetilde{k} = \max \{3, k\}$.
We here summarize the main properties of the space $\VhE$
(we refer to \cite{BM:2012,projectors} for a deeper analysis).

\begin{itemize}
\item [\textbf{(P1)}] \textbf{Polynomial inclusion:} $\Pk_k(E) \subseteq \VhE$;

\item [\textbf{(P2)}] \textbf{Degrees of freedom:}
the following linear operators $\mathbf{D_V}$ constitute a set of DoFs for $\VhE$:
\begin{itemize}
\item[$\mathbf{D_V1}$] the value of $v(\boldsymbol{\xi})$ at any vertex $\boldsymbol{\xi}$ of the polygon $E$,
\item[$\mathbf{D_V2}$] the value of $h_{\boldsymbol{\xi}} \partial_{x_1}v(\boldsymbol{\xi})$ and $h_{\boldsymbol{\xi}} \partial_{x_2}v(\boldsymbol{\xi})$ at any vertex $\boldsymbol{\xi}$ of the polygon $E$,
\item[$\mathbf{D_V3}$] the values of $v$ at $k_e = \max\{0, k-3\}$ distinct points of every edge $e \in \Sigma_h^E$,
\item[$\mathbf{D_V4}$] the values of $h_e \partial_{\nn}v$ at $k_{\nn} = \max \{0, k-2\}$ distinct points of every edge $e \in \Sigma_h^E$,
\item[$\mathbf{D_V5}$] the moments of $v$ against a polynomial basis  $\{m_i\}_i$ of $\Pk_{k-4}(E)$ with $\|m_i\|_{L^{\infty}(E)} = 1$:
$$
\frac{1}{|E|}\int_E v \, m_{i} \, {\rm d}E  \,.
$$
\end{itemize}
Therefore the dimension of $\VhE$ is 
\[
\dim(\VhE) = (3  + k_e + k_{\nn}) n_e +  \max\left\{0, \frac{(k-3)(k-2)}{2} \right\} \,.
\]
\item [\textbf{(P3)}] \textbf{Polynomial projections:}
the DoFs $\mathbf{D_V}$ allow us to compute the following linear operators:
\[
\begin{aligned}
&\Pi^{0,E}_{k-2} \colon D^2 \VhE \to [\Pk_{k-2}(E)]^{2 \times 2}\,, \\
&\Pi^{0,E}_{k-2} \colon \Delta \VhE \to \Pk_{k-2}(E) \,, \\
&\Pi^{0,E}_{k-1} \colon \nabla \VhE \to [\Pk_{k-1}(E)]^2\,, \\
&\Pi^{0,E}_k \colon \VhE \to \Pk_{k}(E) \,.
\end{aligned}
\]
\end{itemize} 

The global  space $\VhO$ is defined by gluing the local spaces with the obvious associated sets of global DoFs:
\begin{equation}
\label{eq:v-glo}
\VhO = \{v \in V \quad \text{s.t.} \quad v_{|E} \in \VhE \quad \text{for all $E \in \Omega_h$} \}\,,
\end{equation}
with dimension
\[
\dim(\VhE) = 3 \, N_V  + (k_e + k_{\nn}) N_e +  \max\left\{0, \frac{(k-3)(k-2)}{2} \right\} N_P \,.
\]

We now recall the optimal approximation properties for the space 
$\VhO$ (see, for instance, \cite{BM:2012,Antonietti-BeiraodaVeiga-Scacchi-Verani:2016}).
\begin{proposition}[Approximation property of $\VhO$]
\label{lm:int VDG}
Under the Assumption \textbf{(A1)} for any $v \in \VhO \cap H^{s}(\Omega_h)$ there exists $v_{\mathcal {I}} \in \VhO$ such that for all $E \in \Omega_h$ it holds
\[
\Vert v - v_{\mathcal {I}} \Vert_{0,E} + 
h_E \vert v - v_{\mathcal {I}} \vert_{1, E} +
h_E^2 \vert v - v_{\mathcal {I}} \vert_{2, E}
\lesssim h_E^{s} |v|_{s,E} \,,
\]
where $2 < s \leq k+1$.
\end{proposition}

\subsection{Virtual Element forms}
\label{sub:forms}

The next step in the construction of our method is to define a discrete version of the continuous forms in \eqref{eq:forms} and \eqref{eq:formar}.
It is clear that for an arbitrary functions in $\VhO$ the forms are not computable since the discrete functions are not known in closed form.
Therefore, following the usual procedure in the VEM setting, we need to construct discrete forms that are computable by the DoFs.

In the light of property \textbf{(P3)} for any $v_h$, $w_h \in \VhE$ we define the computable local discrete bilinear forms:
\begin{equation}
\label{eq:forms-VEM}
\begin{aligned}
& a_h^{D^2,E}(v_h, w_h) = \int_E (\Pi_{k-2}^{0,E} D^2 v_h) : (\Pi_{k-2}^{0,E} D^2 w_h) \, {\rm d}E + h_E^{-2} \mathcal{S}^E(v_h, w_h) \,, 
\\
& a_h^{0,E}(v_h, w_h) = \int_E (\Pi_{k}^{0,E} v_h) (\Pi_{k}^{0,E}  w_h) \, {\rm d}E +  h_E^2 \mathcal{S}^E(v_h, w_h) \,, 
\\
& b_h^{E}(v_h, w_h) = \int_E {\bf u}  \cdot (\Pi_{k-1}^{0,E}\nabla v_h) \ (\Pi_{k}^{0,E}w_h) \, {\rm d}E
\,,
\end{aligned}
\end{equation}
and for any $v_h$, $w_h$, $z_h \in \VhE$ the semi-linear forms 
\begin{equation}
\label{eq:formar-VEM}
\begin{aligned}
&l_h^E(f;v_h,w_h) = \int_E \lambda (f - \Pi^{0,E}_k v_h)\Pi^{0,E}_k w_h  \, {\rm d}E \,,
\\
&r_h^E(z_h; v_h, w_h) = \int_E \phi'(\Pi^{0,E}_k z_h) (\Pi^{0,E}_{k-1}\nabla v_h) \cdot (\Pi^{0,E}_{k-1}\nabla w_h) \, {\rm d}E \,.
\end{aligned}
\end{equation}
The VEM stabilizing term in \eqref{eq:forms-VEM} is given by
\begin{equation}
\label{eq:stab1}
\mathcal{S}^E(v_h, w_h) = 
S^E \bigl( (I - \Pi^{0,E}_k) v_h, \, (I - \Pi^{0,E}_k) w_h \bigr) \,,
\end{equation}
where $S^E(\cdot, \cdot) \colon \VhE \times \VhE \to \R$ is a computable symmetric discrete form satisfying for all $v_h \in \VhE \cap \ker(\Pi^{0,E}_k)$ the following bounds
\begin{equation}
\label{eq:stab2}
\begin{gathered}
\vert v_h \vert^2_{2,E} \lesssim h_E^{-2}S^E(v_h, v_h) \lesssim \vert v_h \vert^2_{2,E} \,,
\\
\Vert v_h \Vert^2_{0,E} \lesssim h_E^2 S^E(v_h, v_h) \lesssim \Vert v_h \Vert^2_{0,E} \,.
\end{gathered}
\end{equation}
Many examples of such stabilization can be found in the VEM literature 
\cite{volley,autostoppisti,BDR:2017,BM:2012}.
In the present paper we consider the so-called \texttt{dofi-dofi} stabilization defined as follows: let $\vec{v}_h$ and $\vec{w}_h$ denote the real valued vectors containing the values of the local degrees of freedom associated to $v_h$, $w_h$ in the space $\VhE$ then
\begin{equation}
\label{eq:dofi}
S^E(v_h, w_h) =\vec{v}_h \cdot \vec{w}_h \,.
\end{equation}
In particular we notice that the linear operators $\mathbf{D_V}$ in property \textbf{(P2)} are properly scaled to recover the bounds \eqref{eq:stab2}.

The global forms can be derived by \eqref{eq:forms-VEM} and \eqref{eq:formar-VEM}, employing the notation in \eqref{eq:XG}.
%
\subsection{Virtual Element problem}
\label{sub:vem problem}

Referring to the space \eqref{eq:v-glo} the discrete bilinear forms  \eqref{eq:forms-VEM} and the discrete semi-linear forms  \eqref{eq:formar-VEM}, we can state the following semi-discrete problems.

\noindent {\bf Advective Cahn-Hilliard VEM problem}: find $c_h(\cdot, t) \in \VhO$ s.t.
\begin{equation}\label{vem-contpbl:1}
\left\{
\begin{aligned}
& a_h^0(\partial_t c_h,v_h) + \frac{ \gamma^2}{\Pe} a_h^{D^2}(c_h,v_h) + \frac{1}{\Pe}r_h(c_h;c_h,v_h) + b_h(c_h,v_h)= 0 \quad \forall v_h \in \VhO ,   \\
& c_h(\cdot,0)=c_{0,h}  \,.
\end{aligned}
\right.
\end{equation}

\noindent {\bf Cahn-Hilliard inpainting VEM problem}:  find $c_h(\cdot, t) \in \VhO$ s.t.
\begin{equation}\label{vem-contpbl:2}
\left\{
\begin{aligned}
& a_h^0(\partial_t c_h,v_h) + { \gamma} a_h^{D^2}(c_h,v_h) + \frac{1}{\gamma} r_h(c_h;c_h,v_h) + l_h(f; c_h, v_h)= 0 \quad \forall v_h \in \VhO , 
\\
& c_h(\cdot,0)=c_{0,h}  \,.
\end{aligned}
\right.
\end{equation}

\noindent
In problems \eqref{vem-contpbl:1} and \eqref{vem-contpbl:2} the discrete initial datum $c_{0,h} \in \VhO$ is the DoFs interpolant of $c_0$, i.e. $\mathbf{D_V}(c_{0,h} - c_0) = \boldsymbol{0}$. \\

In the next step we formulate a fully discrete version of problems \eqref{vem-contpbl:1} and \eqref{vem-contpbl:2}.
We introduce a sequence of time steps $t_n = n \tau$, $n=0,\dots,N$, with time step size $\tau$. 
Next, we define $v_{h, \tau}^n \approx v_h(\cdot, t_n)$ as the
approximation of the function $v_h(\cdot, t) \in \VhO$ at time~$t_n$, $n=0,\dots,N$.
Here we chose  the backward Euler method.
The fully discrete systems consequently reads as follows.

\medskip

\noindent {\bf Advective Cahn-Hilliard discrete problem}.
\begin{equation}\label{fully-contpbl:1}
\left\{
\begin{aligned}
& \text{given $c_{h, \tau}^0 = c_{0,h}$, find $c_{h, \tau}^n$ with $n=1, \dots, N$ s.t. $\forall v_h \in \VhO$ it holds:}
\\
& \frac{1}{\tau}a_h^0(c_{h, \tau}^n - c_{h, \tau}^{n-1},v_h) + 
\frac{ \gamma^2}{\Pe} a_h^{D^2}(c_{h, \tau}^n,v_h) + \frac{1}{\Pe}r_h(c_{h, \tau}^n;c_{h, \tau}^n,v_h) + b_h(c_{h, \tau}^n,v_h)= 0  \,.
\end{aligned}
\right.
\end{equation}

\noindent {\bf Cahn-Hilliard inpainting discrete problem}.
\begin{equation}\label{fully-contpbl:2}
\left\{
\begin{aligned}
& \text{given $c_{h, \tau}^0 = c_{0,h}$, find $c_{h, \tau}^n$ with $n=1, \dots, N$ s.t. $\forall v_h \in \VhO$ it holds:}
\\
& \frac{1}{\tau}a_h^0(c_{h, \tau}^n - c_{h, \tau}^{n-1},v_h) + 
{ \gamma} a_h^{D^2}(c_{h, \tau}^n,v_h) 
+ \frac{1}{\gamma} r_h(c_{h, \tau}^n;c_{h, \tau}^n,v_h) + 
l_h(f; c_{h, \tau}^n, v_h)= 0    \,.
\end{aligned}
\right.
\end{equation}

\section{Numerical results}
\label{sec:num}
In this section, we numerically explore the efficacy of the conforming virtual element discretizations \eqref{fully-contpbl:1} and \eqref{fully-contpbl:2}. In particular, the results of the approximation of the advective Cahn-Hilliard problem  are reported in Section~\ref{S:ACH}, while the ones obtained with the Cahn-Hilliard inpainting problem are collected in Section~\ref{S:CHI}. \\

We remark that the resulting nonlinear systems (\ref{fully-contpbl:1}) and (\ref{fully-contpbl:2}) at each time step are solved by the Newton method, using the $l^2$-norm of the relative residual as a stopping criterion,
with tolerance $\texttt{1e-6}$. Except otherwise stated, the Jacobian linear system is solved by GMRES,
preconditioned by a Block-Jacobi preconditioner,
using the $l^2$-norm of the relative residual as a stopping criterion,
with tolerance $\texttt{1e-8}$.\\

\begin{figure}[ht]
\begin{center}
\subfigure[Cartesian ({\bf QUAD}) mesh]{\includegraphics[width=3.8cm]{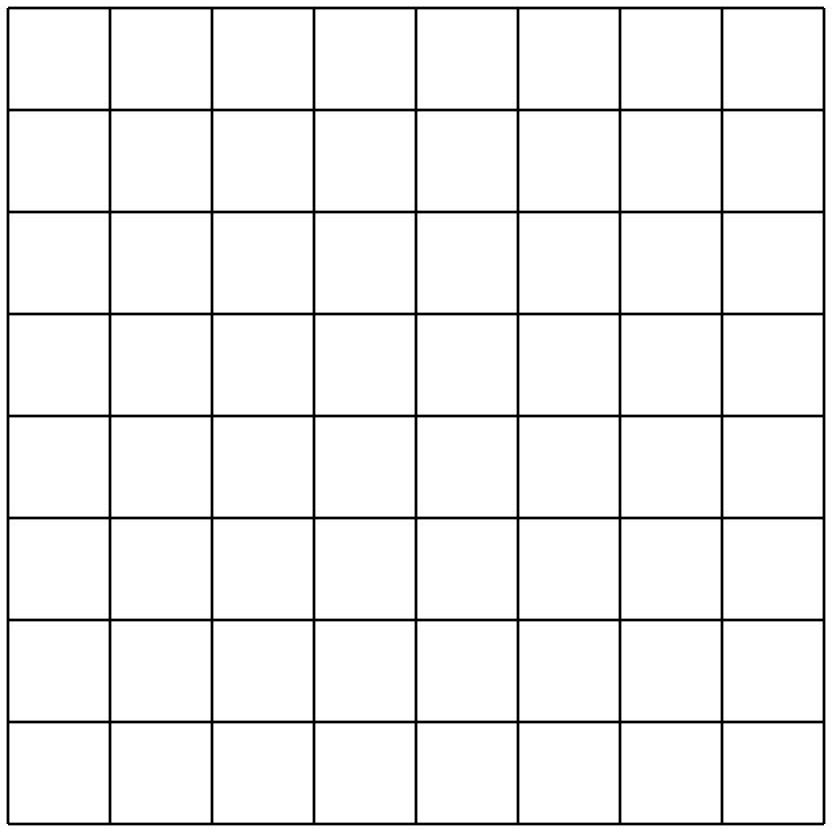}} 
\subfigure[Triangular ({\bf TRI}) mesh]{\includegraphics[width=3.8cm]{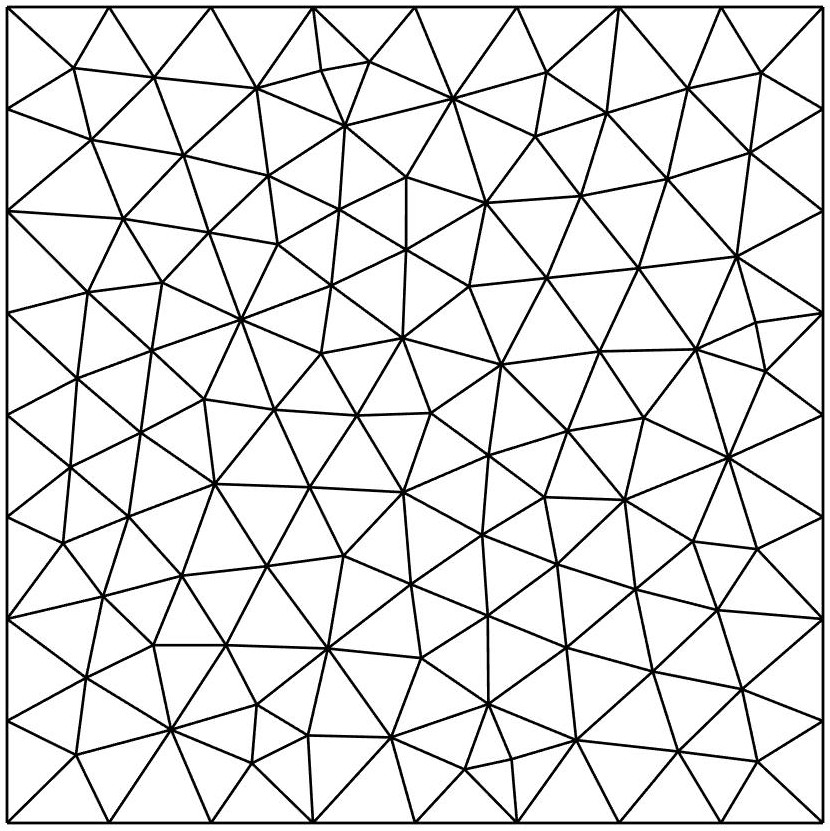}} 
\subfigure[CVT ({\bf CVT}) mesh]{\includegraphics[width=3.8cm]{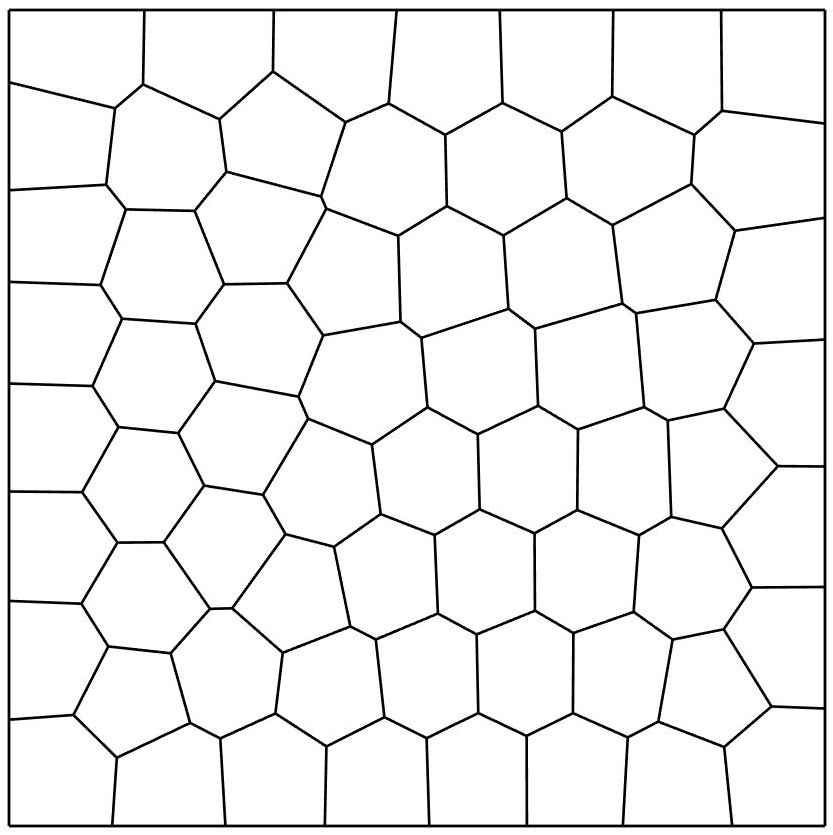}} 
\caption{Example of polygonal meshes used in the numerical tests.}
\label{fig_mesh}
\end{center}
\end{figure}

\begin{table}[h]
\begin{center}
\begin{tabular}{c|c|c|c|c}
\hline
mesh		& $1/h$		& \# elements		& \# nodes	& \# DoFs \\
\hline
{\bf QUAD}	& 128		& 16384			& 16641		& 49923 \\
{\bf TRI}	& 128		& 56932			& 28723		& 86169 \\
{\bf CVT}	& 128		& 16384			& 32943		& 98829 \\
\hline
\end{tabular}
\caption{Mesh size parameter $h$, number of elements, number of nodes and number of degrees of freedom (DoFs) of the polygonal meshes used in the numerical tests.}
\label{table_mesh}
\end{center}
\end{table} 

For the computational mesh, we consider three different mesh families, i.e., quadrilateral ({\bf QUAD}),
triangular ({\bf TRI}) and central Voronoi tessellation ({\bf CVT}) meshes. An example
of a mesh of each family is shown in Figure~\ref{fig_mesh}. The corresponding number of elements, number of nodes, and number of degrees of freedom of the  meshes used in all tests (except Test \ref{parallelTest}) are reported in Table~\ref{table_mesh}.\\

Finally, the simulations have been performed using an in-house Fortran90 parallel code based on the PETSc library \cite{petsc-user-ref}.
Except otherwise stated, the parallel tests
were run on 32 cores of the INDACO linux cluster at the University of Milan (indaco.unimi.it).

\subsection{Advective Cahn-Hilliard problem}
\label{S:ACH}
We consider two scenarios: the evolution of a cross (Tests 1 and 2, Figure~\ref{cross_Pe_100_gam_1d100}) and 
a spinoidal decomposition (Test 3, Figure~\ref{spin2_Pe_100_gam_1d100}).  
In both cases, the convective field $\mathbf{u}$ is taken from \cite{Kay-Styles-Suli:2009}, i.e.
\[
\mathbf{u}(x,y)=f(r)(2y-1,1-2x)^T, \quad (x,y)\in \Omega=(0,1)^2
\]
where
\[
f(r)=\frac{1}{2}\left(1+\tanh\left(\beta\left(\frac{1}{2}-\epsilon-r\right)\right)\right) \quad \mbox{ and } \quad
r=\sqrt{\left(x-\frac{1}{2}\right)^2+\left(y-\frac{1}{2}\right)^2},
\]
with $\beta=200$ and $\epsilon=0.1$.
The parameters $\Pe$ and $\gamma$ in system (\ref{eq:CHforte:1}) are set to 100 and 0.01, respectively.
\begin{table}[!htb]
\begin{center}
\begin{tabular}{c|cc|ccc|ccc|ccc}
\hline
\multicolumn{12}{c}{{\bf Advective Cahn-Hilliard problem, evolution of a cross}} \\
\hline
\hline
\multicolumn{12}{c}{{\bf QUAD mesh with 147456 elements, DoFs = 444675}} \\
\hline
$p$ &\multicolumn{2}{c|}{Mumps}         &\multicolumn{3}{c|}{BJ}                &\multicolumn{3}{c}{GAMG}               &\multicolumn{3}{c}{bAMG}\\
    & nit       & $T_{sol}$      & nit   & it    &$ T_{sol}$    & nit   & it    & $T_{sol}$    & nit   & it    & $T_{sol}$  \\
\hline
1   & 2.2       & 67.8           & 2.2   & 14.7  & 15.5         & 2.2   & 11.5  & 27.8         & 2.2   & 13.3  & 36.7 \\
2   & 2.2       & 39.4           & 2.2   & 33.9  & 10.9         & 2.2   & 13.8  & 19.4         & 2.2   & 13.4  & 21.1 \\
4   & 2.2       & 24.4           & 2.2   & 37.4  & 8.6          & 2.2   & 13.8  & 10.1         & 2.2   & 13.8  & 12.9 \\
8   & 2.2       & 21.6           & 2.2   & 45.8  & 11.7         & 2.2   & 14.2  & 12.7         & 2.2   & 14.0  & 13.3 \\
16  & 2.2       & 14.3           & 2.2   & 46.2  & 6.4          & 2.2   & 14.4  & 7.3          & 2.2   & 14.0  & 8.1  \\
32  & 2.2       & 7.8            & 2.2   & 45.0  & 1.1          & 2.2   & 14.4  & 1.7          & 2.2   & 14.1  & 2.3  \\
48  & 2.2       & 7.2            & 2.2   & 43.7  & 0.82         & 2.2   & 14.8  & 1.3          & 2.2   & 14.2  & 1.7  \\
\hline
\hline
\multicolumn{12}{c}{{\bf CVT mesh with 147456 elements, DoFs = 884814}} \\
\hline
$p$ &\multicolumn{2}{c|}{Mumps}         &\multicolumn{3}{c|}{BJ}                &\multicolumn{3}{c}{GAMG}               &\multicolumn{3}{c}{bAMG}\\
    & nit       & $T_{sol}$      & nit   & it    &$ T_{sol}$    & nit   & it    & $T_{sol}$    & nit   & it    & $T_{sol}$  \\
\hline
1	& OoM	& OoM		& 2.2	& 32.2	& 44.9		& 2.2	& 18.9	& 88.7		& 2.2	& 21.1	& 135.7 \\
2	& 2.2	& 202.1		& 2.2	& 79.6	& 36.1		& 2.2	& 25.6	& 67.8		& 2.2	& 24.1	& 172.2 \\
4	& 2.2	& 123.9		& 2.2	& 95.9 	& 27.6		& 2.2	& 28.5	& 59.8		& 2.2	& 25.6	& 125.3 \\
8	& 2.2	& 85.4		& 2.2	& 107.4	& 23.9		& 2.2	& 30.4	& 45.9		& 2.2	& 26.9	& 74.6 \\
16	& 2.2	& 53.2		& 2.2	& 110.6	& 14.6		& 2.2	& 30.7	& 39.6		& 2.2	& 27.2	& 38.1 \\
32	& 2.2	& 32.4		& 2.2	& 109.7	& 4.7		& 2.2	& 31.2	& 34.3		& 2.2	& 27.1	& 18.1 \\
48	& 2.2	& 27.2		& 2.2	& 108.3	& 3.2		& 2.2	& 30.8	& 30.5		& 2.2	& 27.2 	& 12.9 \\
\hline
\hline
\end{tabular}
\caption{Strong scaling test on {\bf QUAD} and {\bf CVT} meshes, Advective Cahn-Hilliard, evolution of a cross. $p$=number of procs; nit=average Newton iterations per time step; it=average GMRES iterations per Newton iteration; $T_{sol}$=average CPU time in seconds per time step; OoM=out of memory.}
\label{table_scal_cvt}
\end{center}
\end{table}

\subsubsection{Test 1: parallel performance of the solver}
\label{parallelTest}

We first study the performance of the parallel solver by comparing four different methods:
\begin{itemize}
\item Mumps: the Jacobian system at each Newton iteration is solved by the parallel direct solver Mumps \cite{amestoy.2001,amestoy.2006};
\item BJ: the Jacobian system at each Newton iteration is solved by the Block-Jacobi preconditioner implemented in the PETSc object PCBJACOBI;
\item GAMG: the Jacobian system at each Newton iteration is solved by the Algebraic Multigrid preconditioner implemented in the PETSc object PCGAMG, with default settings;
\item bAMG: the Jacobian system at each Newton iteration is solved by the Algebraic Multigrid preconditioner boomerAMG \cite{henson-yang} of the HYPRE library \cite{hypre}.
\end{itemize}

\begin{figure}[ht]
\begin{center}
{\bf QUAD} mesh\\
\subfigure[$t = 0.01$]{\includegraphics[width=2.5cm]{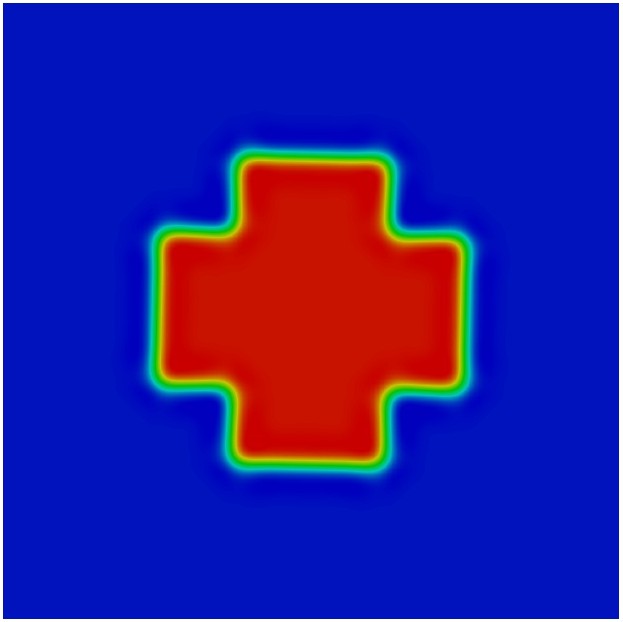}}
\subfigure[t = 0.1]{\includegraphics[width=2.5cm]{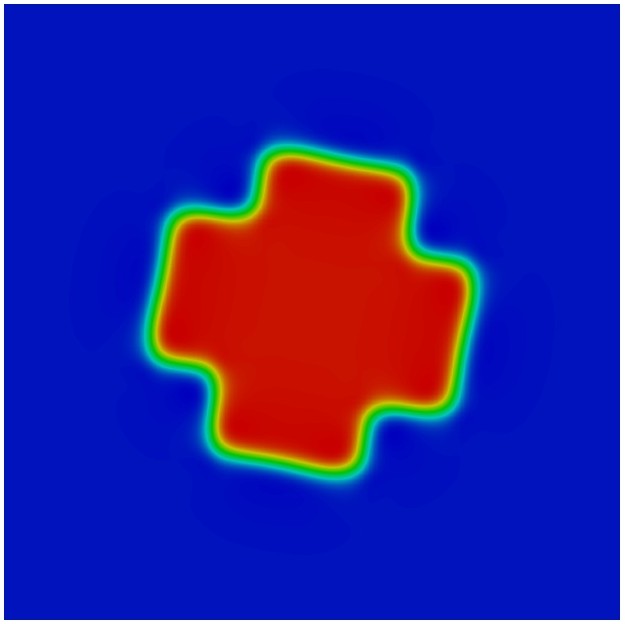}}
\subfigure[$t = 1$]{\includegraphics[width=2.5cm]{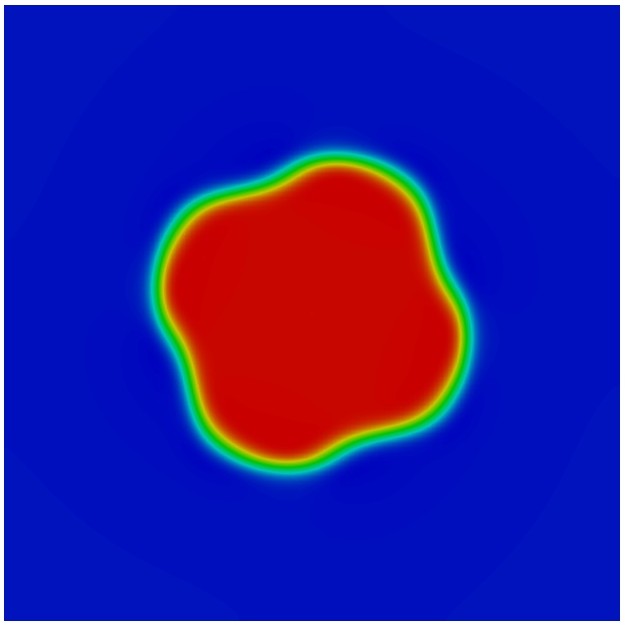}}
\subfigure[$t = 10$]{\includegraphics[width=2.5cm]{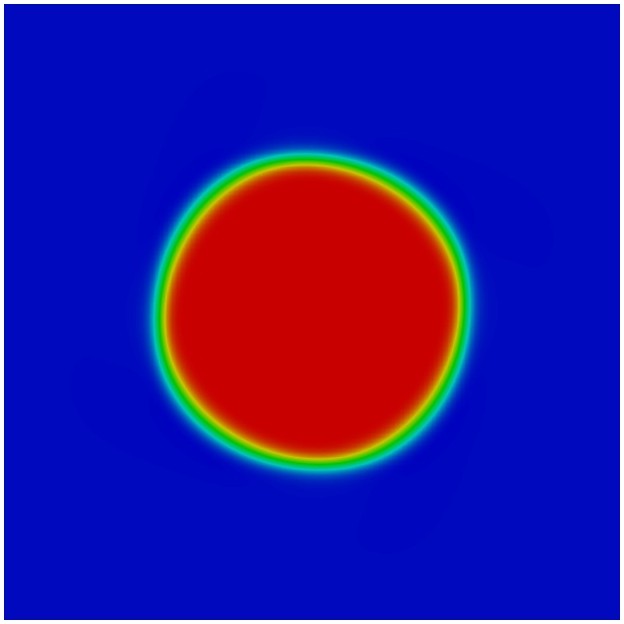}}\\
{\bf TRI} mesh\\
\subfigure[$t = 0.01$]{\includegraphics[width=2.5cm]{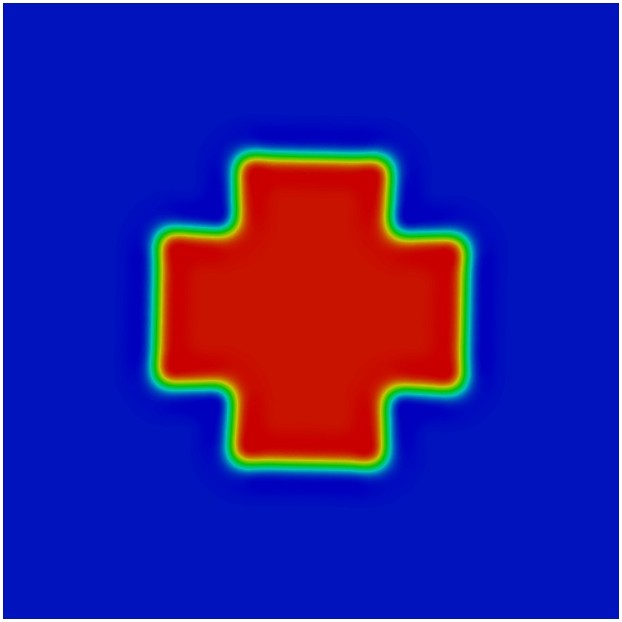}}
\subfigure[t = 0.1]{\includegraphics[width=2.5cm]{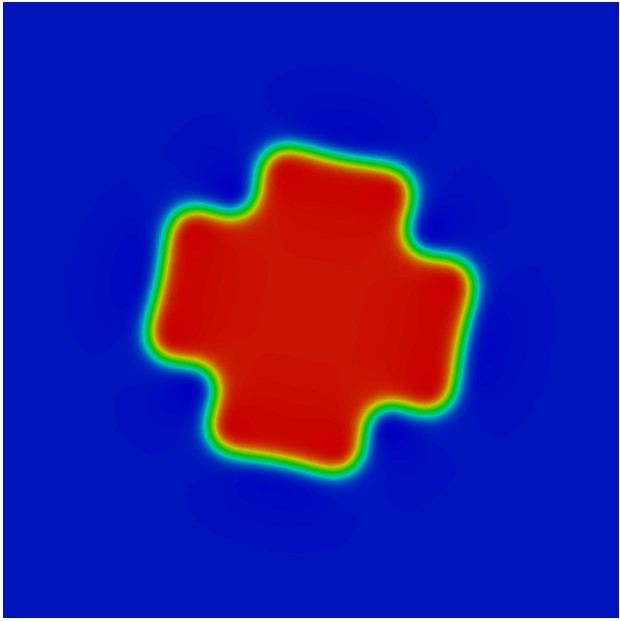}}
\subfigure[$t = 1$]{\includegraphics[width=2.5cm]{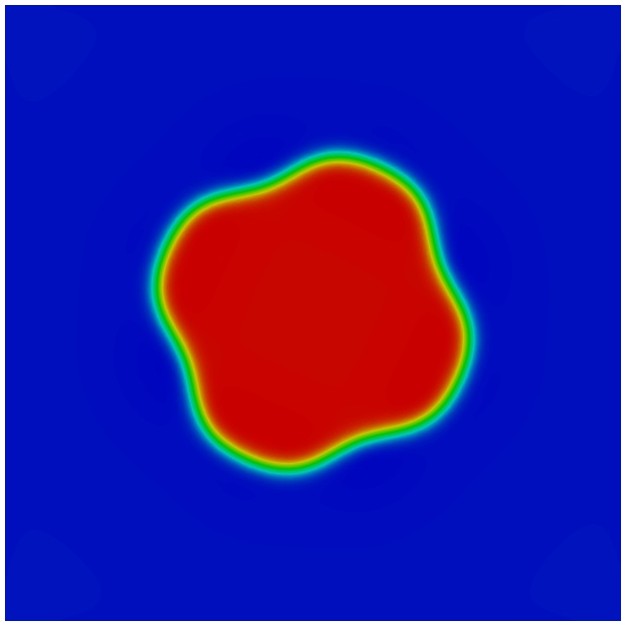}}
\subfigure[$t = 10$]{\includegraphics[width=2.5cm]{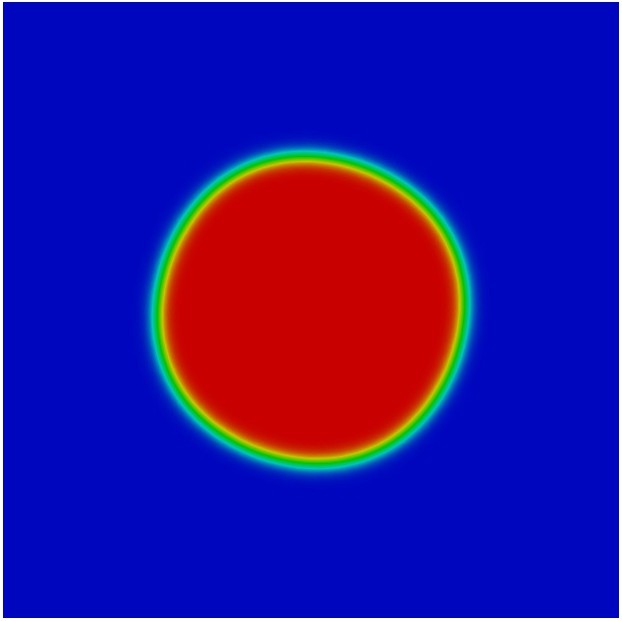}}\\
{\bf CVT} mesh\\
\subfigure[$t = 0.01$]{\includegraphics[width=2.5cm]{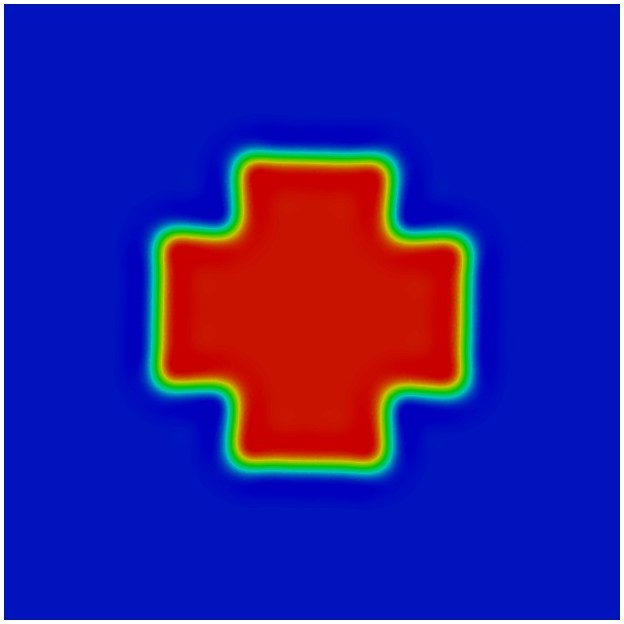}}
\subfigure[t = 0.1]{\includegraphics[width=2.5cm]{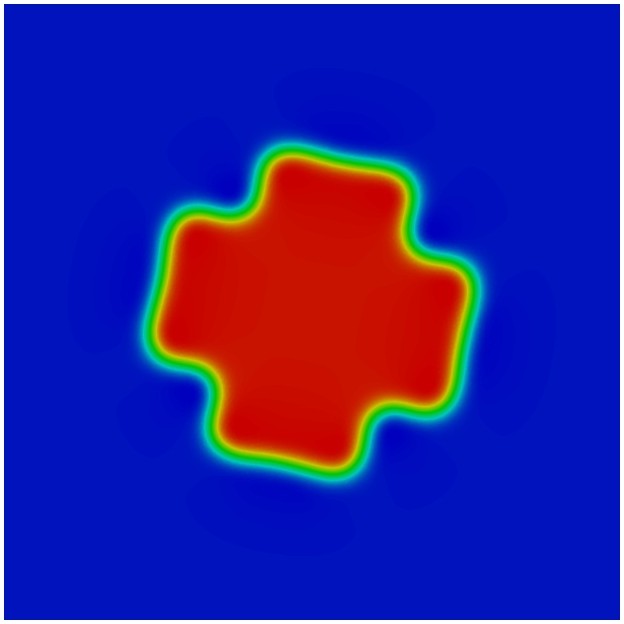}}
\subfigure[$t = 1$]{\includegraphics[width=2.5cm]{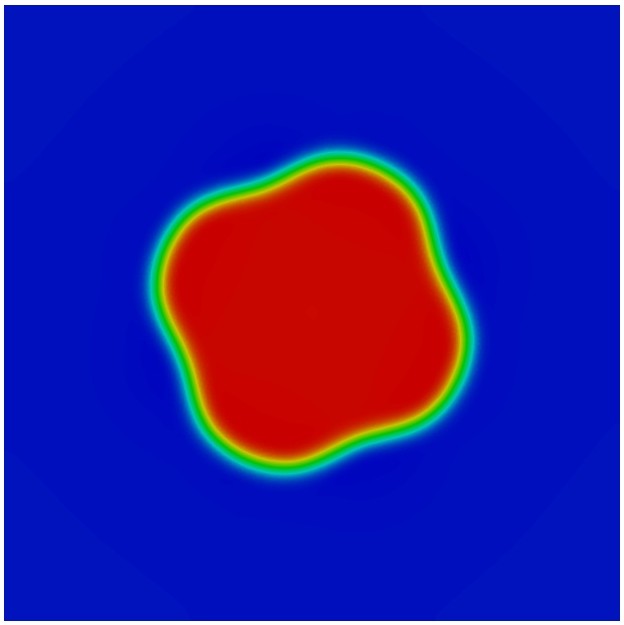}}
\subfigure[$t = 10$]{\includegraphics[width=2.5cm]{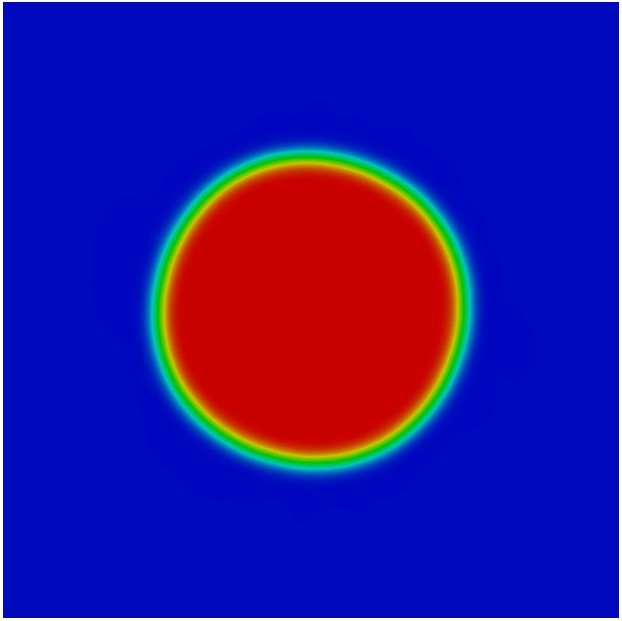}}\\
\caption{Test 2, evolution of a cross with convection on the unit square, $\gamma=1/100$, $\Pe=100$. 
The mesh parameters are reported in Table~\ref{table_mesh}.}
\label{cross_Pe_100_gam_1d100}
\end{center}
\end{figure}

The initial datum $c_0$ is a piecewise constant function whose
jump set has the shape of a cross, see Figure~\ref{cross_Pe_100_gam_1d100}, Panels (a-e-i).
The time step size considered is $\tau=2e-5$ and the simulation is run for 50 time steps, up to $T=1e-3$. 
The unit square domain is discretized by a {\bf QUAD} mesh of 147456 elements ($1/h=384$, DoFs = 444675) 
and a {\bf CVT} mesh of 147456 elements ($1/h=384$, DoFs = 884814); see Table~\ref{table_scal_cvt}.  
We increase the number of processors from 1 to 48, keeping fixed the global number of DoFs, thus performing
a strong scaling test. The code is run on the Galileo100 cluster of CINECA laboratory (www.cineca.it).\\

The results show that the four parallel solvers are all scalable, since the CPU times reduce 
when the number of processors increase. As expected, the Algebraic Multigrid preconditioners exhibit a scalable
behavior of GMRES iterations, which remain almost constant with respect to the number of processors.
The BJ preconditioner shows an initial increase in terms of iterations, but after 8-16 processors they remain stable.
We believe that this scalable behavior of the BJ preconditioner is due to the dominant effect of the mass matrix,
which improves the conditioning of the Jacobian linear system. Indeed, the most effective solver results to be
the BJ preconditioner, which in case of the {\bf CVT} mesh is about 9 times as fast as Mumps, 10 times as fast as GAMG and
4 times as fast as bAMG.

\subsubsection{Test 2: evolution of a cross under convection}

As in the previous test, the initial datum $c_0$ is again a piecewise constant function whose
jump set has the shape of a cross.
The time step size considered is $\tau=2e-5$ and the simulation is run for $500000$ time steps, up to $T=10$.
The evolution of the cross simulated on the three computational meshes with data reported in Table~\ref{table_mesh} 
is displayed in Figure~\ref{cross_Pe_100_gam_1d100}. The cross, rotating under the convective field, evolves towards a circle.

\begin{figure}[ht]
\begin{center}
{\bf QUAD} mesh\\
\subfigure[$t = 0.01$]{\includegraphics[width=2.5cm]{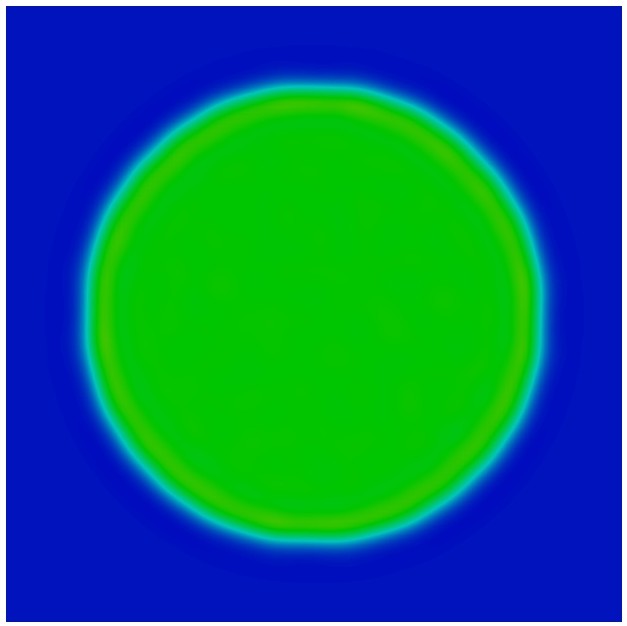}}
\subfigure[$t = 1$]{\includegraphics[width=2.5cm]{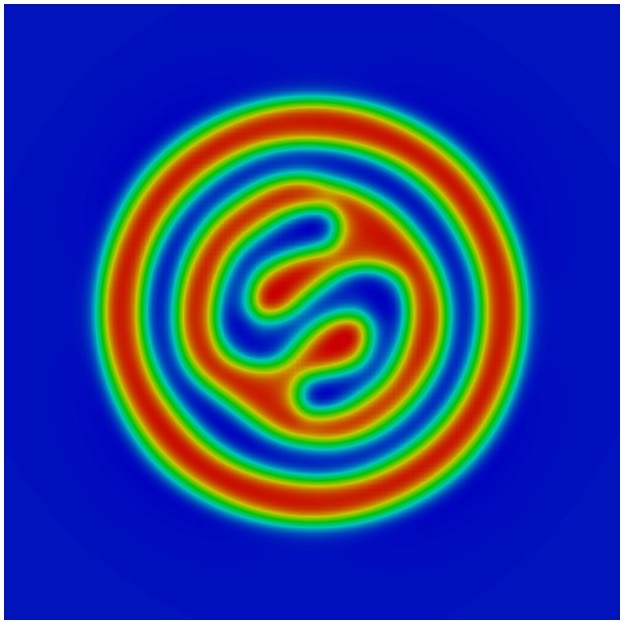}}
\subfigure[$t = 5$]{\includegraphics[width=2.5cm]{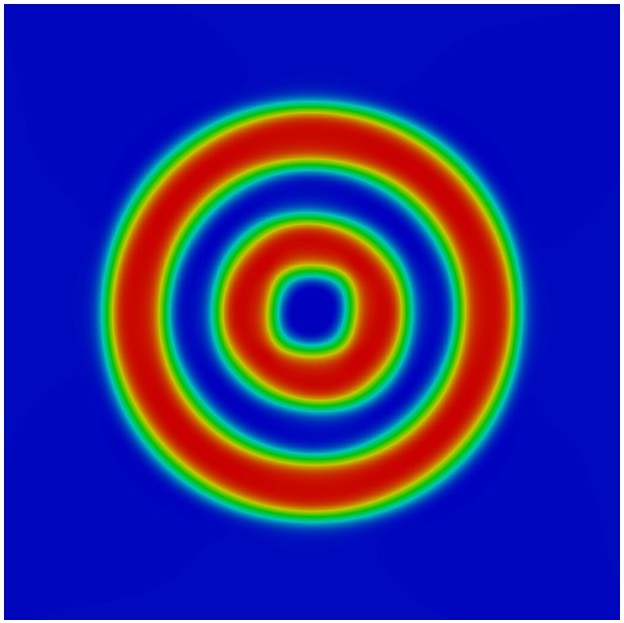}}
\subfigure[$t = 10$]{\includegraphics[width=2.5cm]{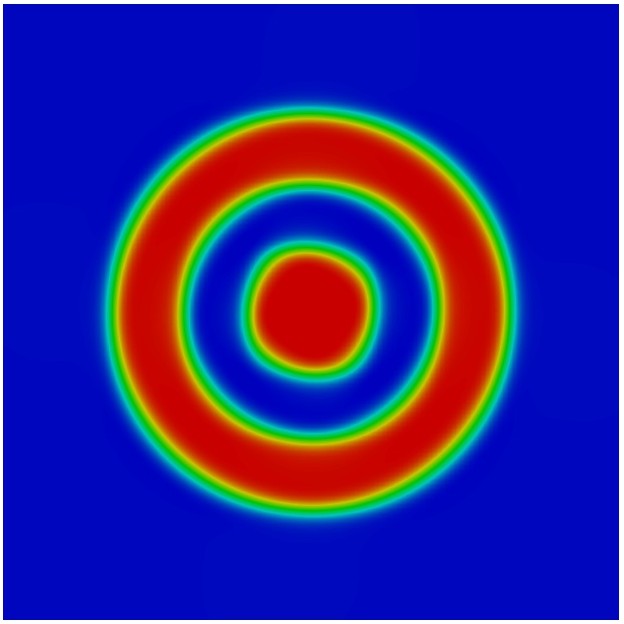}}\\
{\bf TRI} mesh\\
\subfigure[$t = 0.01$]{\includegraphics[width=2.5cm]{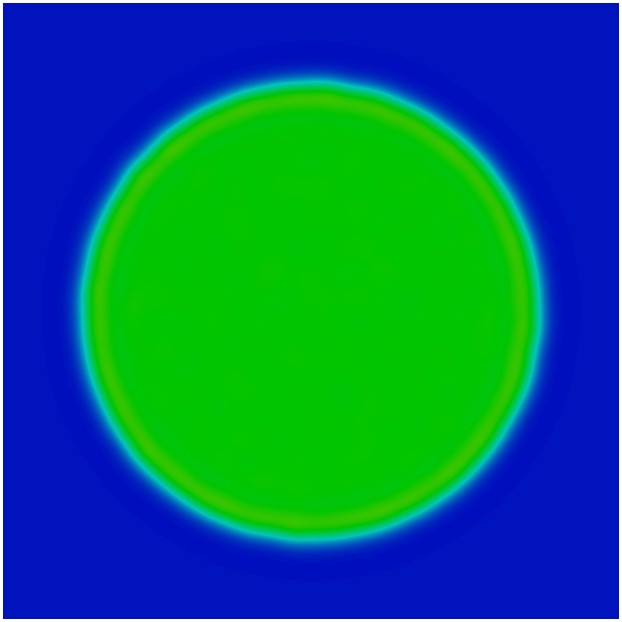}}
\subfigure[$t = 1$]{\includegraphics[width=2.5cm]{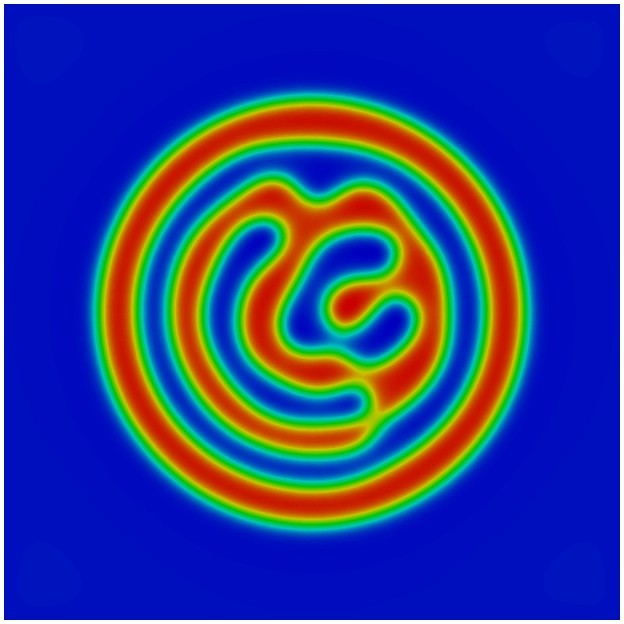}}
\subfigure[$t = 5$]{\includegraphics[width=2.5cm]{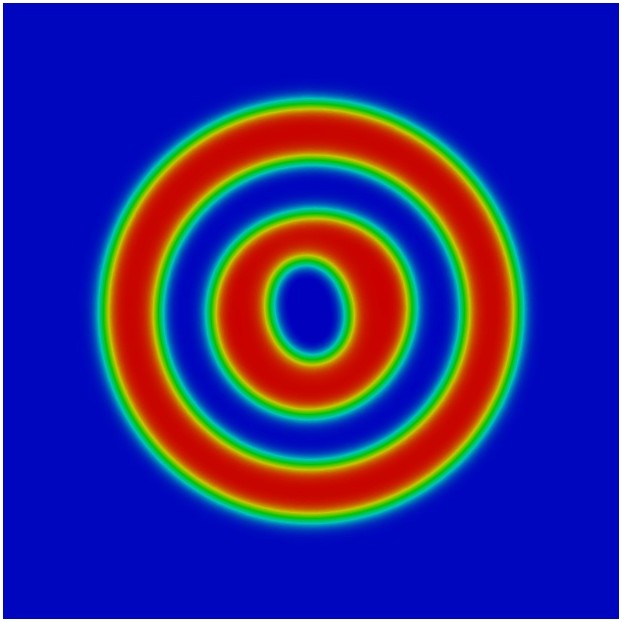}}
\subfigure[$t = 10$]{\includegraphics[width=2.5cm]{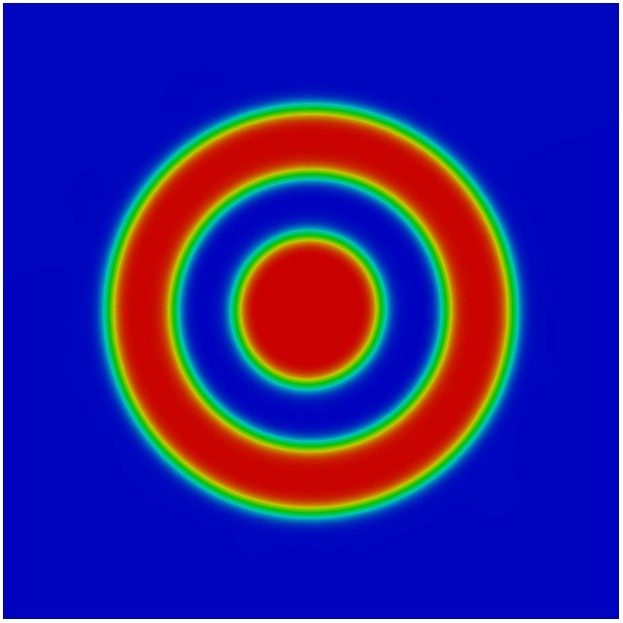}}\\
{\bf CVT} mesh\\
\subfigure[$t = 0.01$]{\includegraphics[width=2.5cm]{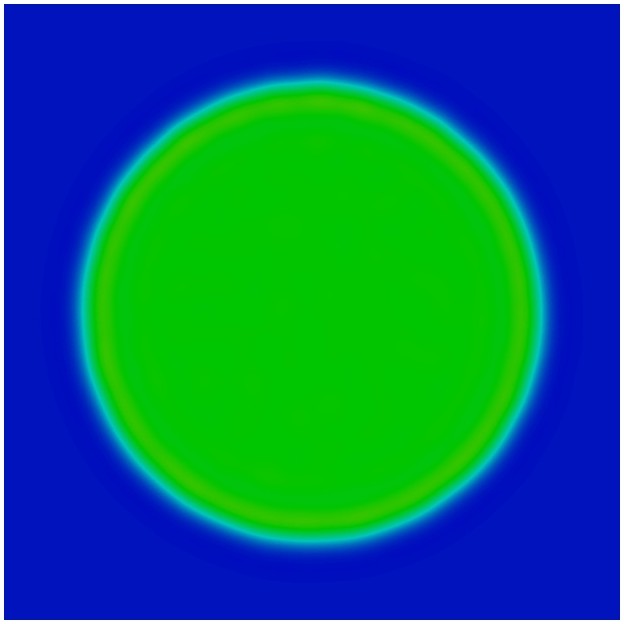}}
\subfigure[$t = 1$]{\includegraphics[width=2.5cm]{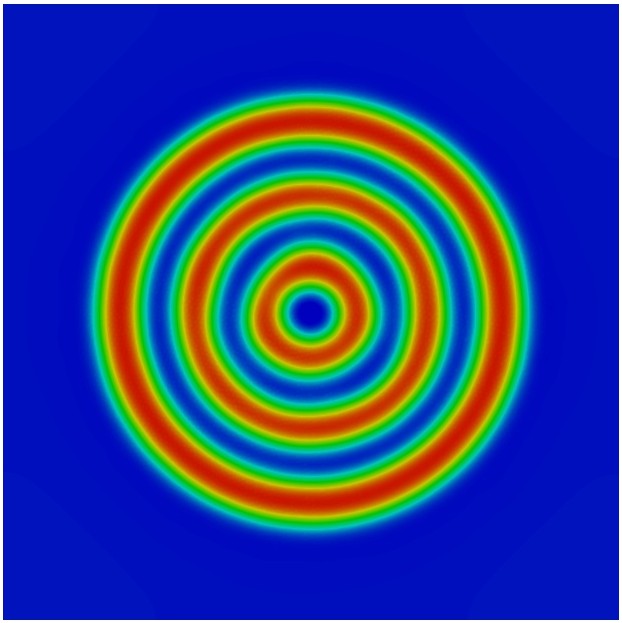}}
\subfigure[$t = 5$]{\includegraphics[width=2.5cm]{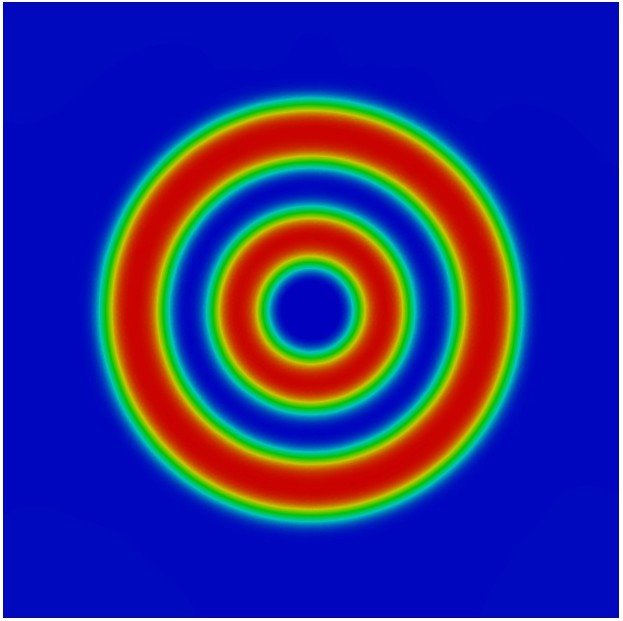}}
\subfigure[$t = 10$]{\includegraphics[width=2.5cm]{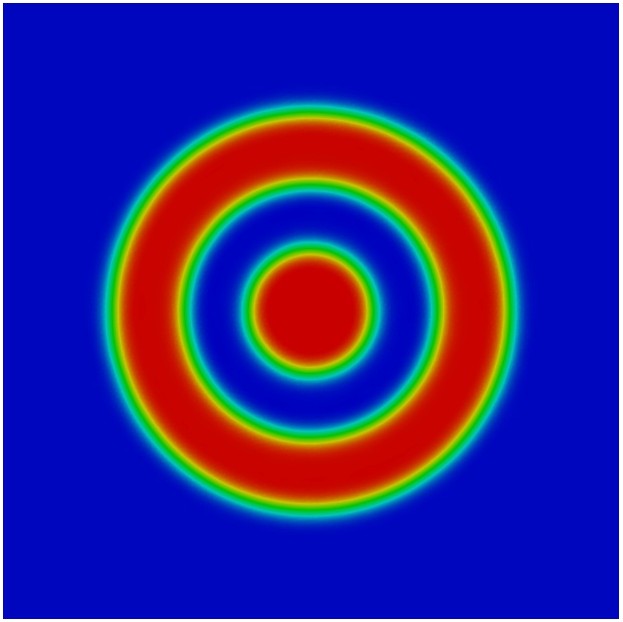}}
\caption{Test 3, spinoidal decomposition of a random disk with convection on the unit square, $\gamma=1/100$, $\Pe=100$. 
The mesh parameters are reported in Table~\ref{table_mesh}.}
\label{spin2_Pe_100_gam_1d100}
\end{center}
\end{figure}

\subsubsection{Test 3: evolution of spinodal decomposition under convection}

The initial datum is now a small uniformly distributed random perturbation about zero, within a circle; 
see Figure~\ref{spin2_Pe_100_gam_1d100}, Panels (a-e-i). 
The time step size considered is $\tau=2e-5$ and the simulation is run for $500000$ time steps, up to $T=10$.
The evolution of the spinoidal decomposition on the three computational meshes 
is displayed in Figure~\ref{spin2_Pe_100_gam_1d100}. 
The initial random distribution evolves very quickly into bulk regions. Then, the convective term makes 
the bulk regions to form concentric circles, which tends very slowly to a central circular bulk region.

\subsection{Cahn-Hilliard inpainting problem}
\label{S:CHI}
\begin{figure}[!htbp]
\begin{center}
{\bf QUAD} mesh\\
\subfigure[$t = 0$]{\includegraphics[width=3cm]{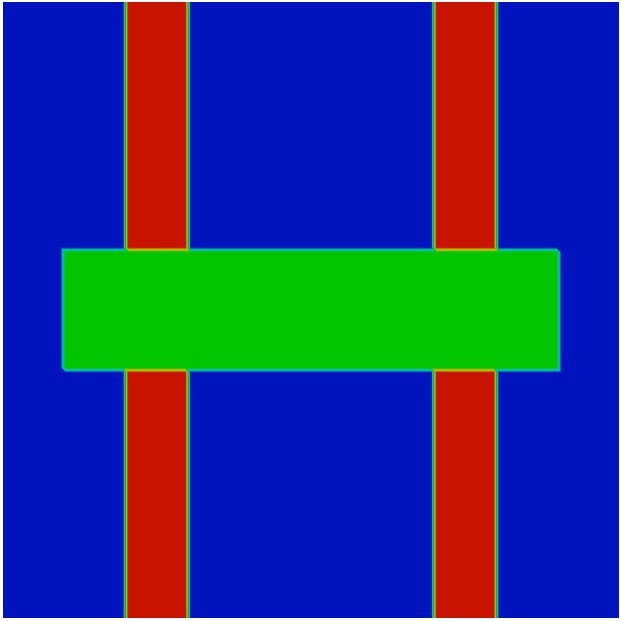}}
\subfigure[$t = 0.02$]{\includegraphics[width=3cm]{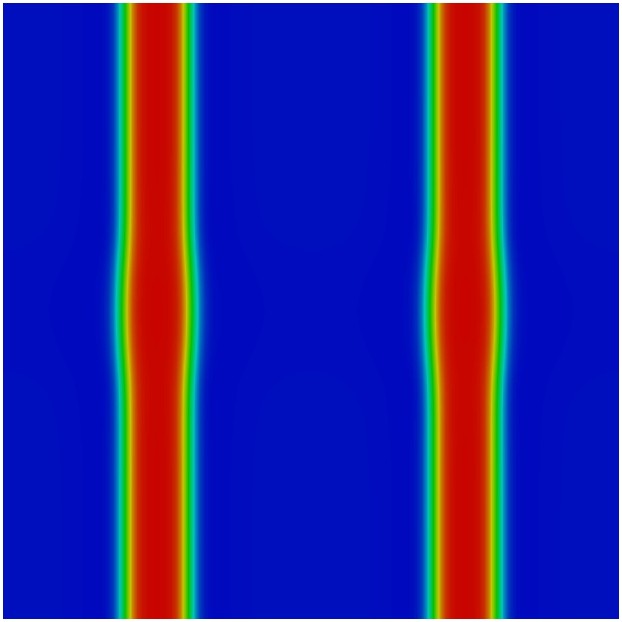}}
\subfigure[$t = 0.02$(binary)]{\includegraphics[width=3cm]{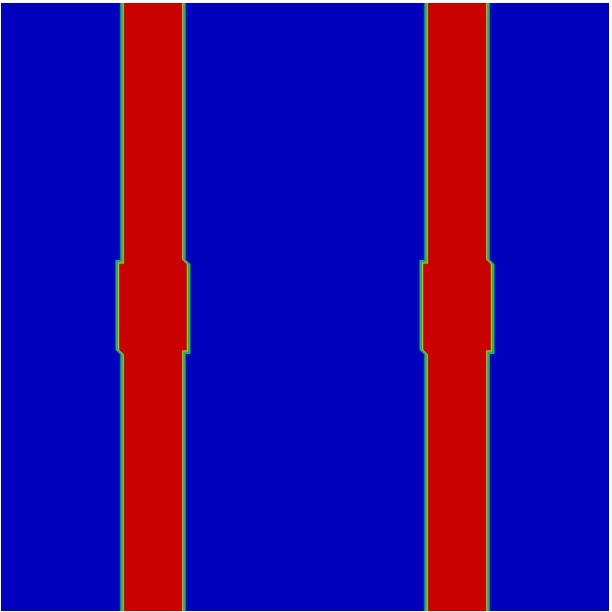}}\\
\medskip
{\bf TRI} mesh\\
\subfigure[$t = 0$]{\includegraphics[width=3cm]{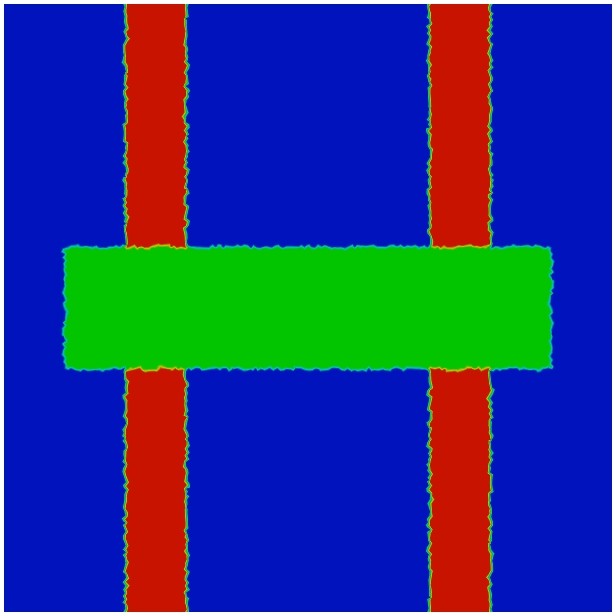}}
\subfigure[$t = 0.02$]{\includegraphics[width=3cm]{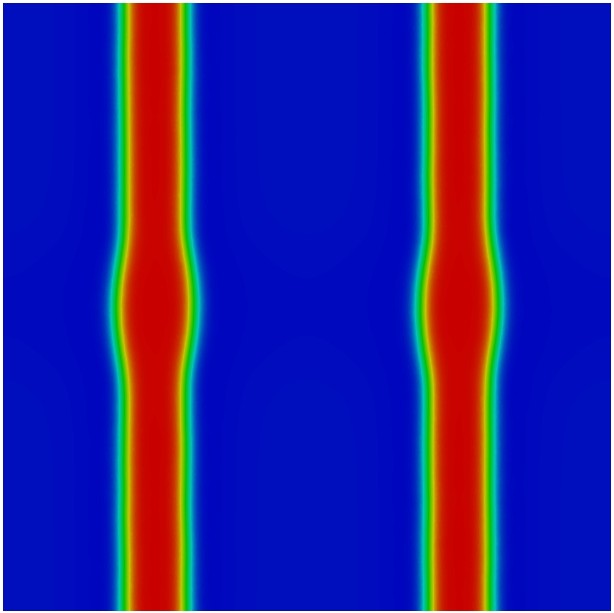}}
\subfigure[$t = 0.02$(binary)]{\includegraphics[width=3cm]{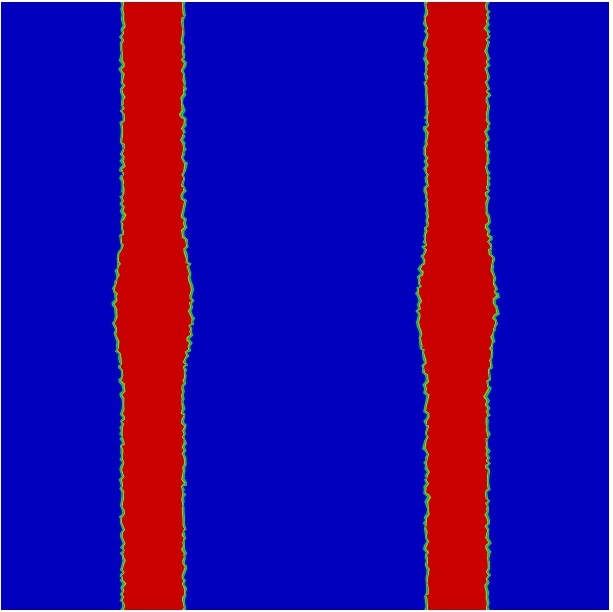}}\\
\medskip
{\bf CVT} mesh\\
\subfigure[$t = 0$]{\includegraphics[width=3cm]{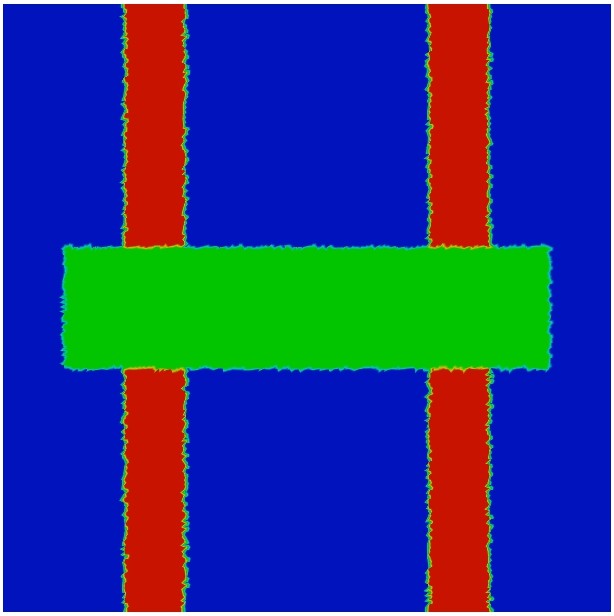}}
\subfigure[$t = 0.02$]{\includegraphics[width=3cm]{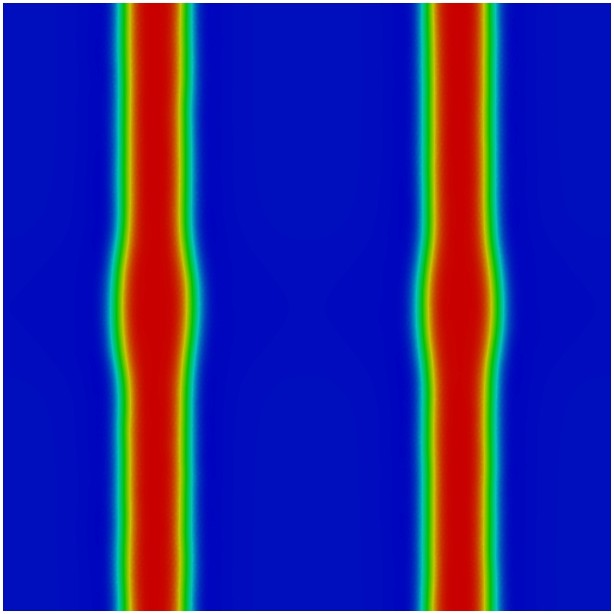}}
\subfigure[$t = 0.02$(binary)]{\includegraphics[width=3cm]{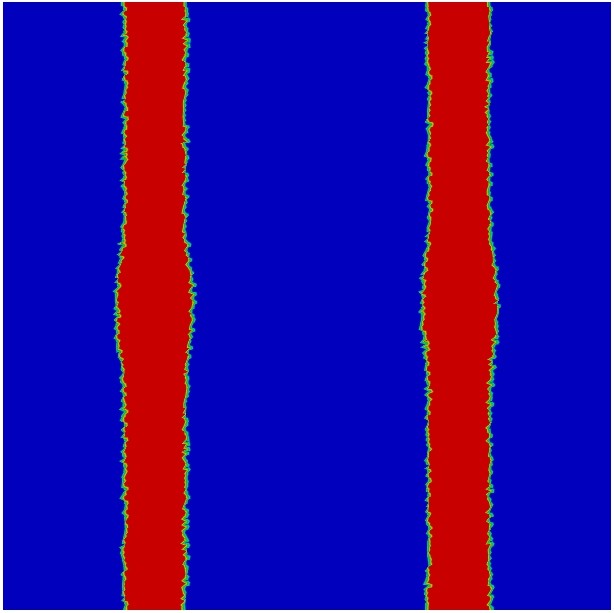}}
\caption{Test 4, impainting of a double stripe. 
The mesh parameters are reported in Table~\ref{table_mesh}.
Left: initial configuration ($t=0$). Middle: final configuration ($t=T=0.02$). Right: final configuration ($t=T=0.02$) without smoothing
effects, projecting the solution $c$ to 0.95 if $c>0$ and to $-0.95$ if $c<0$.}
\label{imp3_fig}
\end{center}
\end{figure}

\begin{figure}[ht]
\begin{center}
{\bf QUAD} mesh\\
\subfigure[$t = 0$]{\includegraphics[width=3cm]{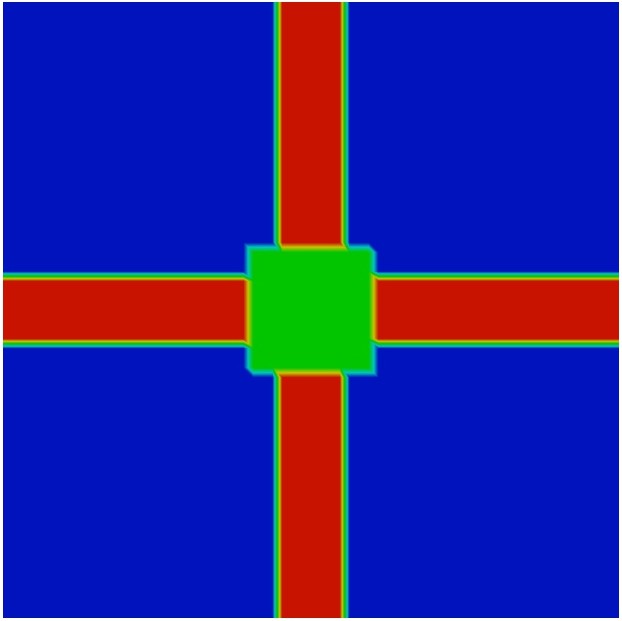}}
\subfigure[$t = 0.02$]{\includegraphics[width=3cm]{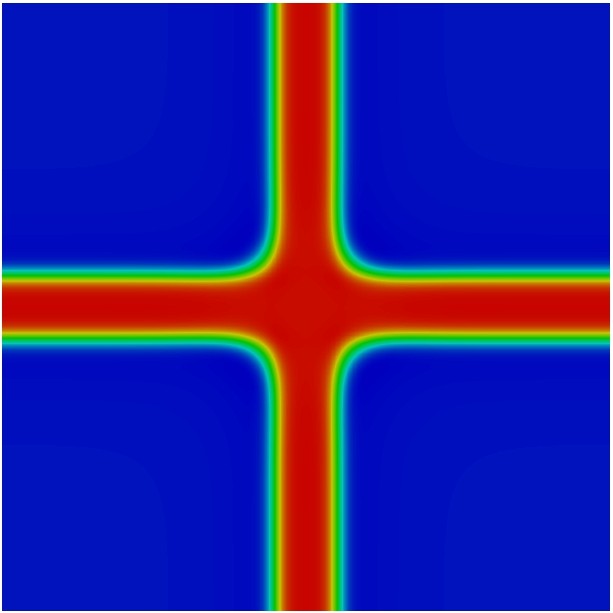}}
\subfigure[$t = 0.02$(binary)]{\includegraphics[width=3cm]{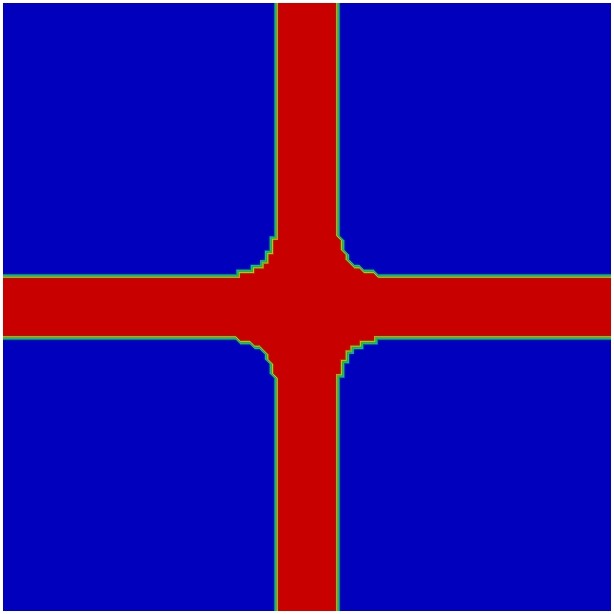}}\\
\medskip
{\bf TRI} mesh\\
\subfigure[$t = 0$]{\includegraphics[width=3cm]{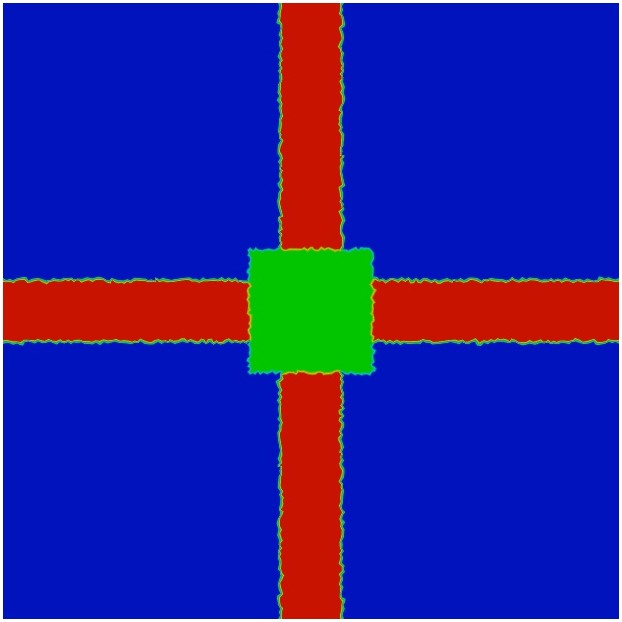}}
\subfigure[$t = 0.02$]{\includegraphics[width=3cm]{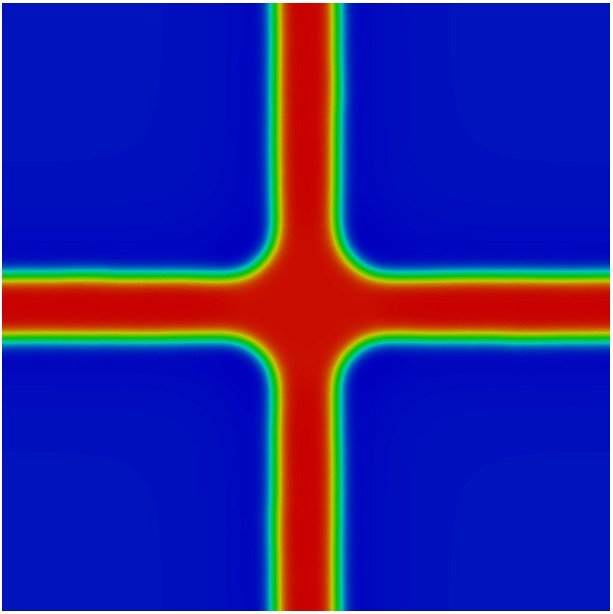}}
\subfigure[$t = 0.02$(binary)]{\includegraphics[width=3cm]{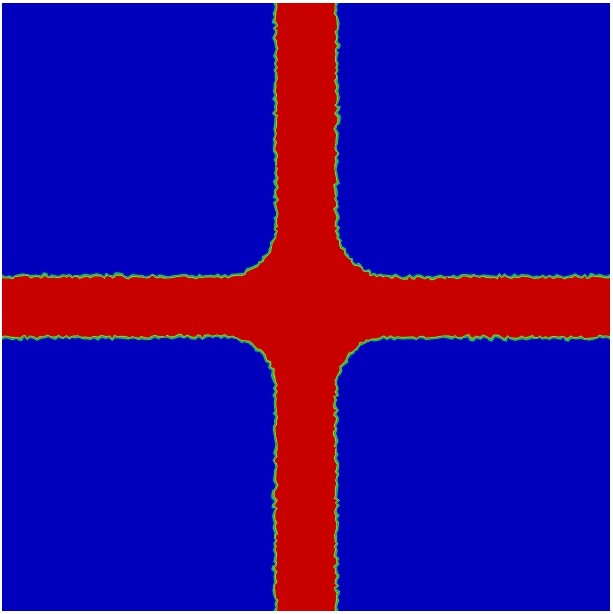}}\\
\medskip
{\bf CVT} mesh\\
\subfigure[$t = 0$]{\includegraphics[width=3cm]{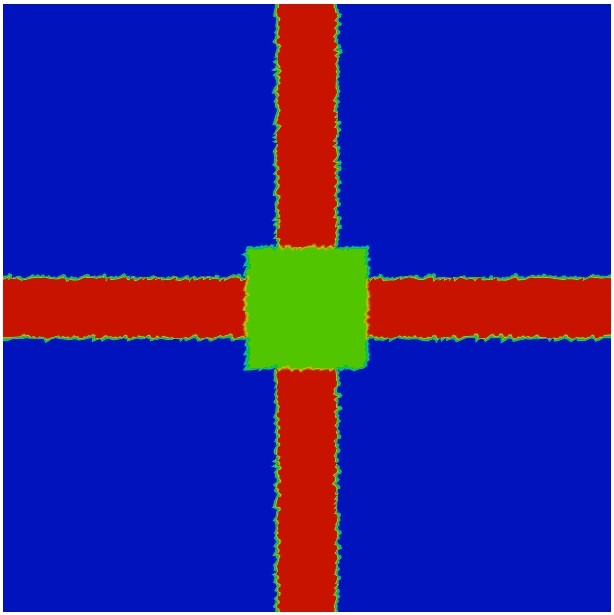}}
\subfigure[$t = 0.02$]{\includegraphics[width=3cm]{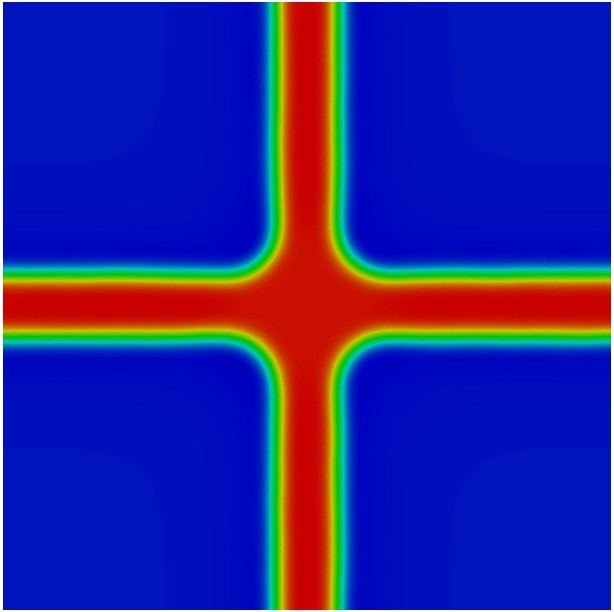}}
\subfigure[$t = 0.02$(binary)]{\includegraphics[width=3cm]{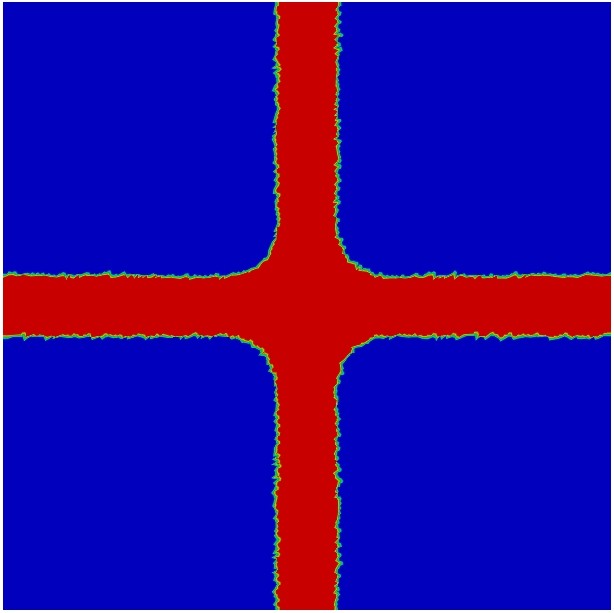}}\\
\caption{Test 5, impainting of a cross.
The mesh parameters are reported in Table~\ref{table_mesh}. 
Left: initial configuration ($t=0$). Middle: final configuration ($t=T=0.02$). Right: final configuration ($t=T=0.02$) without smoothing
effects, projecting the solution $c$ to 0.95 if $c>0$ and to $-0.95$ if $c<0$.}
\label{imp2_fig}
\end{center}
\end{figure}

We consider three scenarios: inpainting of a double stripe (Test 4, Figure~\ref{imp3_fig}), 
inpaiting of a cross (Test 5, Figure~\ref{imp2_fig}) and inpaiting of a circle (Test 6, Figure~\ref{imp5_fig}).
The time step size considered is $\tau=2e-5$ and the simulation is run for 50 time steps, up to $T=1e-3$.
The parameters $\gamma$ and $\lambda_0$ are set to 0.01 and 50000, respectively.
In all next tests, the time step size considered is $\tau=2e-5$ and the simulation is run for $1000$ time steps, up to $T=0.02$.

\subsubsection{Test 4: inpainting of a double stripe}
In this test, the initial configuration consists of two vertical stripes with a central horizontal damage, see Figure~\ref{imp3_fig}.
At the final instant $t=T=0.02$, the correct double stripe configuration is recovered, for all mesh configurations. 
We show also the final configuration without smoothing effects, projecting the solution $c$ to 0.95 if $c>0$ and to $-0.95$ if $c<0$ (binary).

\subsubsection{Test 5: inpainting of a cross}
Here, the initial configuration consists of two stripes, one vertical and one horizontal, crossing at the center of the domain, with a central square damage, see Figure~\ref{imp2_fig}. At the final instant $T=0.02$, the correct cross configuration is recovered.  As before, we also report  the final configuration without smoothing effects, projecting the solution $c$ to 0.95 if $c>0$ and to $-0.95$ if $c<0$ (binary).

\begin{figure}[!htbp]
\begin{center}
{\bf QUAD} mesh\\\
\subfigure[$t = 0$]{\includegraphics[width=3cm]{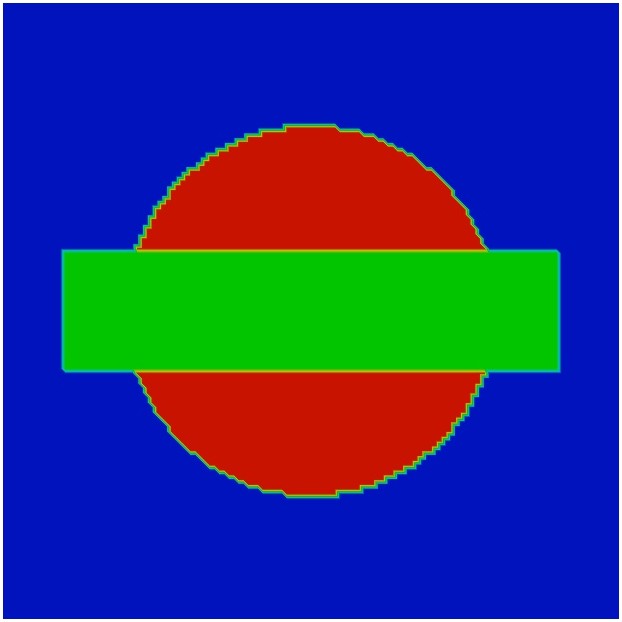}}
\subfigure[$t = 0.02$]{\includegraphics[width=3cm]{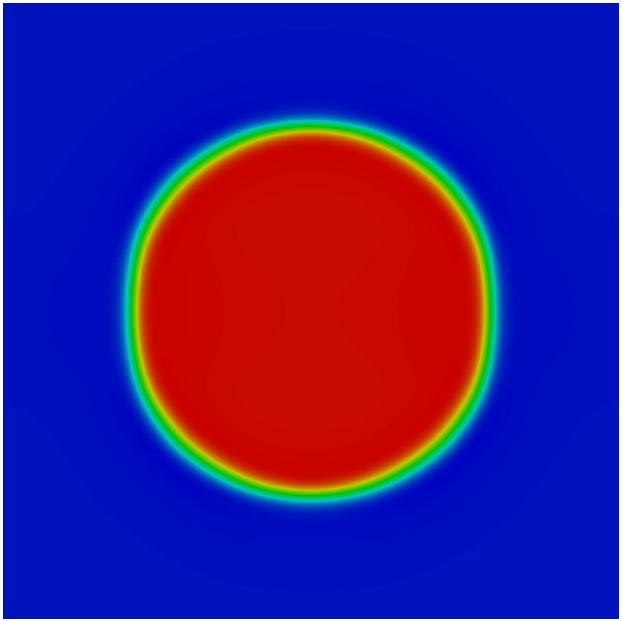}}
\subfigure[$t = 0.02$ (binary)]{\includegraphics[width=3cm]{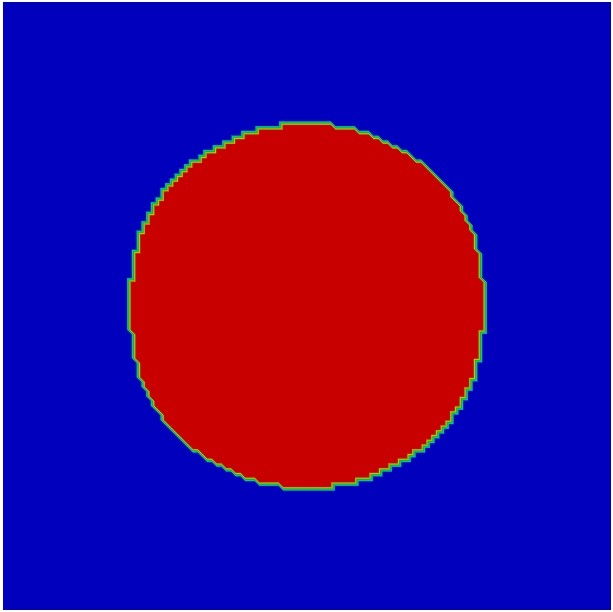}}\\
\medskip
{\bf TRI} mesh\\
\subfigure[$t = 0$]{\includegraphics[width=3cm]{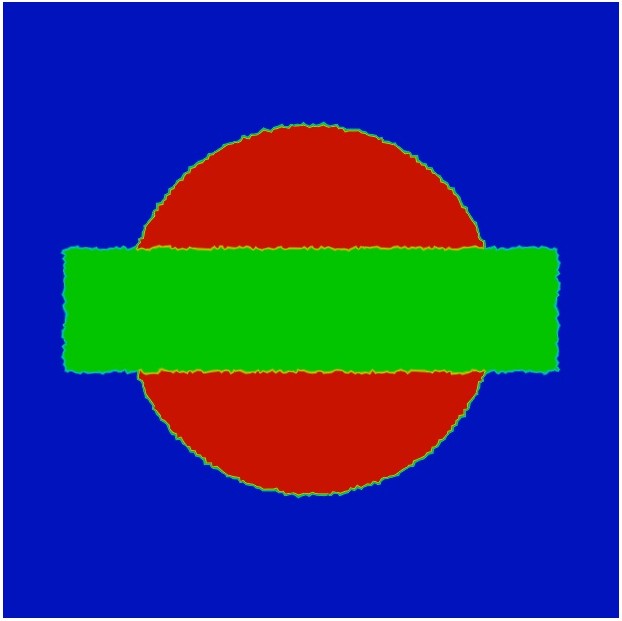}}
\subfigure[$t = 0.02$]{\includegraphics[width=3cm]{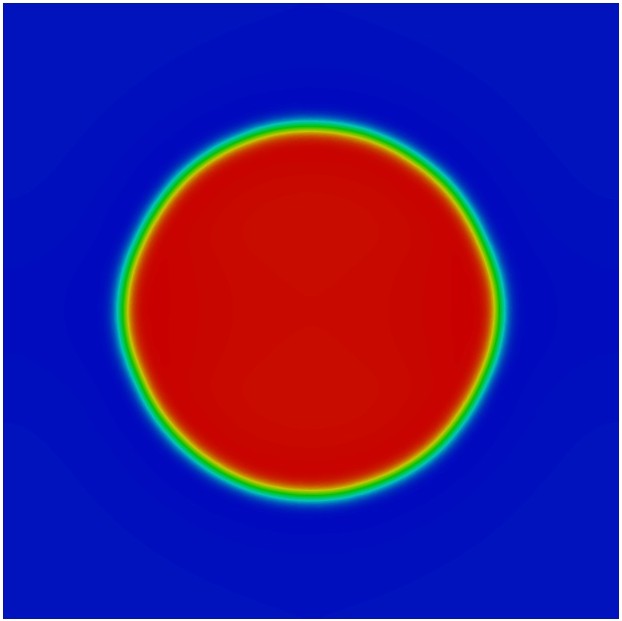}}
\subfigure[$t = 0.02$(binary)]{\includegraphics[width=3cm]{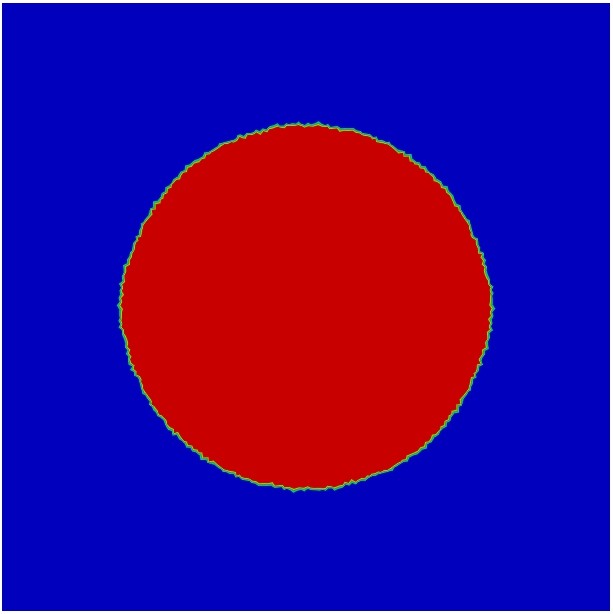}}\\
\medskip
{\bf CVT} mesh\\
\subfigure[$t = 0$]{\includegraphics[width=3cm]{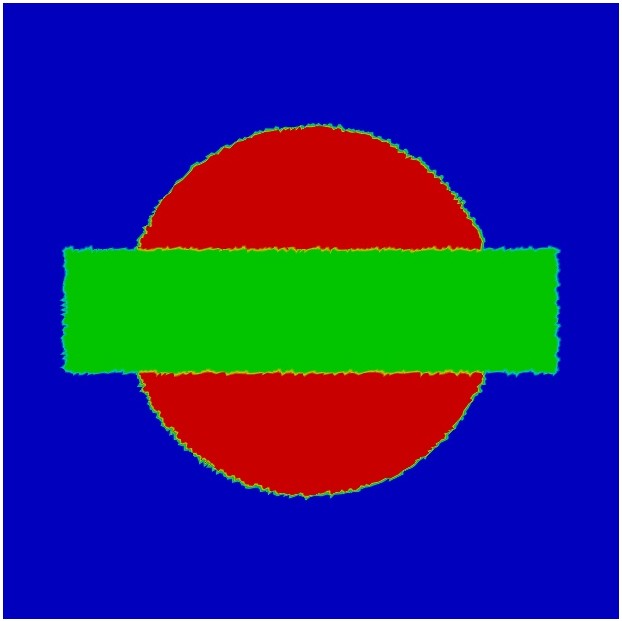}}
\subfigure[$t = 0.02$]{\includegraphics[width=3cm]{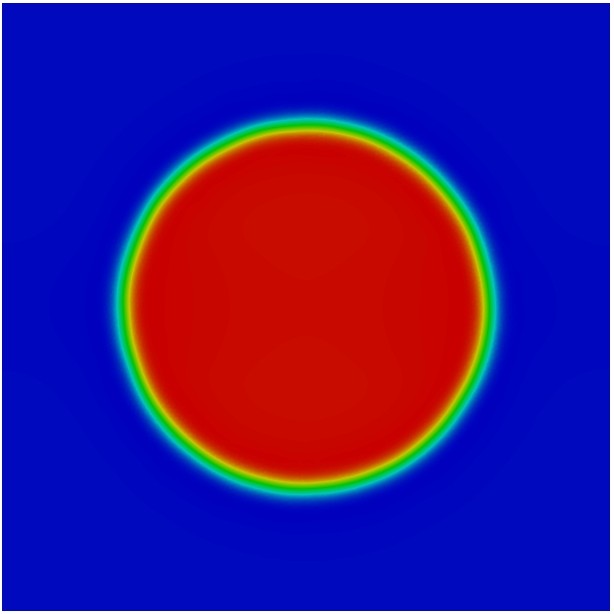}}
\subfigure[$t = 0.02$(binary)]{\includegraphics[width=3cm]{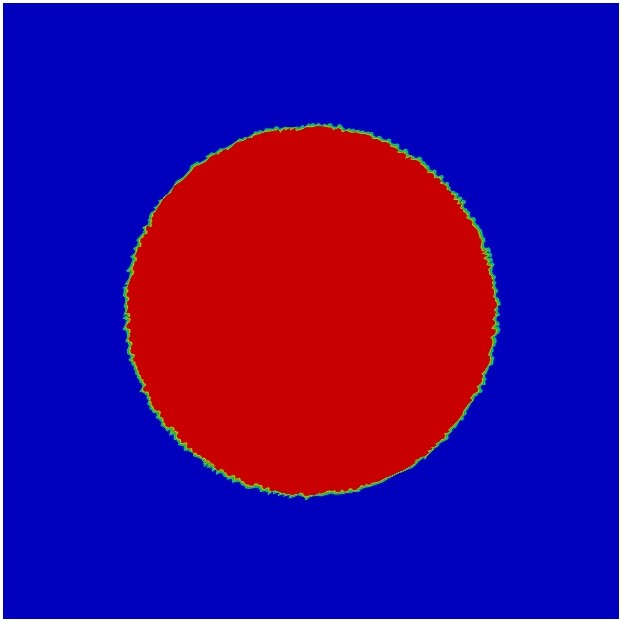}}
\caption{Test 6, impainting of a circle. 
The mesh parameters are reported in Table~\ref{table_mesh}. Left: initial configuration ($t=0$). Middle: final configuration ($t=T=0.02$). Right: final configuration ($t=T=0.02$) without smoothing effects, projecting the solution $c$ to 0.95 if $c>0$ and to $-0.95$ if $c<0$.}
\label{imp5_fig}
\end{center}
\end{figure}
\subsubsection{Test 6: inpainting of a circle}
In the final test, the initial configuration is a circle with a horizontal central damage, see Figure~\ref{imp5_fig}.
At the final instant $T=0.02$, the correct circle configuration is recovered, for all mesh configurations. 
This can be appreciated also from Figure~\ref{imp5_fig} (right panel) where we report  the final configuration (binary plot) projecting the solution $c$ to 0.95 if $c>0$ and to $-0.95$ if $c<0$.

\FloatBarrier
\section{Conclusions}\label{S:conclusions}
In this paper we considered the $C^1$-Virtual Element conforming approximation on polygonal meshes of some variants of the Cahn-Hilliard equation. In particular, we focused on the advective Cahn-Hilliard problem and the Cahn-Hilliard impainting problem. In the first part of the paper we introduced the continuous problems and we gave a detailed description of the virtual element discretizations, while in the second part we numerically explored the efficacy of the proposed methodology through a wide campaign of numerical experiments.

\section*{Acknowledgements} 
P.F. Antonietti and M. Verani have been partially funded by MIUR PRIN research grants n. 201744KLJL and n. 20204LN5N5. The Authors are members of INdAM  GNCS.
The Authors are also grateful to the INDACO scientific platform of the University of Milan and to the CINECA laboratory for the usage of the Galileo100 cluster.

\bibliographystyle{plain}
\bibliography{VEM}

\end{document}